\documentclass[journal, onecolumn, 12pt]{IEEEtran}

\usepackage{graphicx}
\usepackage{amssymb, algorithm, algpseudocode}
\usepackage{amsmath, color, epsf}

\newtheorem{proposition}{Proposition}

\begin{document}

\title{Parallel Algorithms for Constrained Tensor Factorization via Alternating Direction Method of Multipliers}

\author{Athanasios~P. ~Liavas,~\IEEEmembership{Member,~IEEE,}
and~Nicholas~D.~Sidiropoulos,~\IEEEmembership{Fellow,~IEEE}
\thanks{Original manuscript received Aug. 30, 2014; revised \today. Part of this work has been accepted for
presentation at {\em IEEE ICASSP 2015} \cite{LiaSid:ICASSP2015}.}
\thanks{A. P. Liavas is with the Department of Electronic and Computer Engineering,
Technical University of Crete, Chania 73100, Greece (email: {\tt liavas@telecom.tuc.gr}).}
\thanks{N. D. Sidiropoulos is with the Department of Electrical and Computer Engineering, University of Minnesota, Minneapolis,
MN 55455, USA (email: {\tt nikos@ece.umn.edu}). His work was supported in part by NSF IIS-1247632.}
}

\markboth{Submitted to the IEEE Transactions on Signal Processing, 2nd Revision} {}

\maketitle

\begin{abstract}
Tensor factorization has proven useful in a wide range of applications, from sensor array processing to communications, speech and audio
signal processing, and machine learning. With few recent exceptions, all tensor factorization algorithms were originally developed for
centralized, in-memory computation on a single machine; and the few that break away from this mold do not easily incorporate practically
important constraints, such as non-negativity. A new constrained tensor factorization framework is proposed in this paper, building upon the
Alternating Direction Method of Multipliers (ADMoM). It is shown that this simplifies computations, bypassing the need to solve
constrained optimization problems in each iteration; and it naturally leads to distributed algorithms suitable for parallel implementation.
This opens the door for many emerging big data-enabled applications.
The methodology is exemplified using non-negativity as a baseline constraint, but the proposed framework can
incorporate many other types of constraints. Numerical experiments are encouraging, indicating that ADMoM-based non-negative
tensor factorization (NTF) has high potential as an alternative to state-of-the-art approaches.
\end{abstract}

\section{Introduction}

Tensor factorization\footnote{In the literature, the terms {\em factorization} and {\em decomposition} are often used
interchangeably, even though the latter alludes to exact decomposition, whereas the former may include
a residual term.} has proven useful in a wide range of signal processing applications,
such as direction of arrival estimation \cite{SidBroGia00}, communication signal intelligence \cite{SidGiaBro00},
and speech and audio signal separation \cite{NioMokSidPot08,Fevotte2010}, as well as cross-disciplinary areas,
such as community detection in social networks \cite{PapFalSid2012},
and chemical signal analysis \cite{BroSidChem1998}. More recently, there has been significant activity in
applying tensor factorization theory and methods to problems in machine learning research - see \cite{psychovgaltis}.

There are two basic tensor factorization models: {\em parallel factor analysis} (PARAFAC) \cite{Har70,Har72} 
also known as {\em canonical decomposition} (CANDECOMP) \cite{CarCha70}, or CP (and CPD) for CANDECOMP-PARAFAC 
(Decomposition), or {\em canonical polyadic decomposition} (CPD, again); and the Tucker3 model 
\cite{Tucker1966}. Both are sum-of-outer-products models which historically served as cornerstones for further 
developments, e.g., block term decomposition \cite{Lathauwer08decompositionsof}, and upon which the vast 
majority of tensor applications have been built. In this paper, we will primarily focus on the CP model.

Whereas for low-enough\footnote{E.g., relative to the sum of Kruskal-ranks of the latent factor matrices. 
Looser bounds can be guaranteed almost-surely.} rank CP is already unique `on its own,' any side information
can (and should) be used to enhance identifiability and estimation performance in practice. Towards
this end, we may exploit known properties of the sought latent factors, such as non-negativity, sparsity,
monotonicity, or unimodality \cite{SmiBroGel}. Whereas many of these properties can be handled with
existing tensor factorization software, they generally complicate and slow down model fitting.

Unconstrained tensor factorization is already a hard non-convex (multi-linear) problem; even rank-one
least-squares tensor approximation is NP-hard \cite{tensorNP}. Many tensor factorization algorithms rely
on alternating optimization, usually alternating least-squares (ALS), and imposing e.g., non-negativity
and/or sparsity entails replacing linear least-squares conditional updates of the factor matrices with 
non-negative and/or sparse least-squares updates. In addition to ALS, many derivative-based methods have 
been developed that update all model parameters at once, see \cite{tomasi2006comparison} and references therein, 
and \cite{Acar2013,Sorber2013} for recent work in this direction.

With few recent exceptions, all tensor factorization algorithms were originally developed for
centralized, in-memory computation on a single machine. This model of computation is inadequate
for emerging big data-enabled applications, where the tensors to be analyzed cannot be loaded
on a single machine, the data is more likely to reside in cloud storage, and cloud computing,
or some other kind of high performance parallel architecture, must be used for the actual computation.

A carefully optimized Hadoop/MapReduce \cite{hadoop,dean2008mapreduce} implementation of the
basic ALS CP-decomposition algorithm was developed in \cite{kang2012gigatensor}, which reported
$100$-fold scaling improvements relative to the prior art. The jist of \cite{kang2012gigatensor}
is to avoid the explicit computation of `blown-up' intermediate matrix products in the ALS algorithm,
particularly for sparse tensors, and parallelization is achieved by splitting the computation of outer
products. On the other hand, \cite{kang2012gigatensor} is not designed for high performance computing
(e.g., mesh) architectures, and it does not incorporate constraints on the factor matrices.

A random sampling approach was later proposed in \cite{PapFalSid2012}, motivated by recent progress in
randomized algorithms for matrix algebra. The idea of \cite{PapFalSid2012} is to create and analyze multiple
randomly sub-sampled parts of the tensor, then combine the results using a common piece of data to anchor the
constituent decompositions. The downside of \cite{PapFalSid2012} is that it only works for sparse tensors, and
it offers no identifiability guarantees - although it usually works well for sparse tensors.

A different approach based on generalized random sampling was recently proposed in
\cite{PARCOMP-ICASSP2014,PARCOMP-SPM2014}. The idea is to create multiple randomly compressed mixtures (instead
of sub-sampled parts) of the original tensor, analyze them all in parallel, and then combine the results. The main
advantages of \cite{PARCOMP-ICASSP2014,PARCOMP-SPM2014} over \cite{PapFalSid2012} are that i) identifiability
can be guaranteed, ii) no sparsity is needed, and iii) there are theoretical scalability guarantees.

Distributed CP decomposition based on the ALS algorithm has been considered in \cite{AlmeidaCAMSAP2013},
and more recently in \cite{AlmeidaICASSPP2014}, which exploit the inherent parallelism in the matrix
version of the linear least-squares subproblems to split the computation in different ways, assuming an
essentially `flat' architecture for the computing nodes. Regular (e.g., mesh) architectures and constraints
on the latent factors are not considered in \cite{AlmeidaCAMSAP2013,AlmeidaICASSPP2014}.

In this paper, we develop algorithms for constrained tensor factorization based on
Alternating Direction Method of Multipliers (ADMoM). ADMoM has recently attracted renewed
interest \cite{ADMM_Boyd_et_al_NOW}, primarily for solving certain types of convex optimization problems
in a distributed fashion. However, it can also be used to tackle non-convex
problems, such as non-negative matrix factorization \cite{ADMM_Boyd_et_al_NOW}, albeit its convergence
properties are far less understood in this case. We focus on non-negative CP decompositions
as a working problem, due to the importance of the CP
model and non-negativity constraints; but our approach can be generalized
to many other types of constraints on the latent factors, as well as other tensor
factorizations, such as Tucker3, and tensor completion.

The advantages of our approach are as follows. First, during each ADMoM iteration, we avoid the solution of
constrained optimization problems, resulting in considerably smaller computational complexity per iteration
compared to constrained least-squares based algorithms, such as alternating non-negative least-squares (NALS).
Second, our approach leads naturally to distributed algorithms suitable for parallel implementation on regular
high-performance computing (e.g., mesh) architectures. Finally, our approach can easily incorporate many other
types of constraints on the latent factors, such as sparsity.

Numerical experiments are encouraging, indicating that ADMoM-based NTF has significant
potential as an alternative to state-of-the-art approaches.

The rest of the manuscript is structured as follows. In Section \ref{Section_NTF}, we present the NTF problem and
in Section \ref{Section_ADMoM} we present the general ADMoM framework. In Section \ref{Section_ADMoM_NTF}, we develop ADMoM
for NTF, while in Section \ref{Section_ADMoM_LNTF} we develop distributed ADMoM for large NTF. In Section
\ref{Section_Numerical}, we 
test the behavior of the developed schemes with numerical
experiments. Finally, in Section \ref{Section_Conclusion}, we conclude the paper.

\subsection{Notation}

Vectors, matrices, and tensors are denoted by small, capital, and underlined capital bold letters, respectively;
for example, ${\bf x}$, ${\bf X}$, and $\underline{\bf X}$. $\mathbb{R}^{I\times J \times K}_+$ denotes the set
of $(I\times J\times K)$ real non-negative tensors, while $\mathbb{R}^{I\times J}_+$ denotes the set of $(I\times J)$ real
non-negative matrices. $\| \cdot \|_F$ denotes the Frobenius norm of the
tensor or matrix argument, ${\bf A}^\dagger$ denotes the Moore-Penrose
pseudoinverse of matrix ${\bf A}$, and $({\bf A})_+$ denotes the projection of matrix ${\bf A}$ onto the set of
element-wise non-negative matrices. 
The outer product of three
vectors ${\bf a} \in \mathbb{R}^{I \times 1}$, ${\bf b} \in \mathbb{R}^{J \times 1}$, and
${\bf c} \in \mathbb{R}^{K \times 1}$ is the rank-one tensor
${\bf a} \circ {\bf b} \circ {\bf c} \in  \mathbb{R}^{I \times J \times K}$ with elements
$({\bf a} \circ {\bf b} \circ {\bf c})(i,j,k) = {\bf a}(i) {\bf b}(j) {\bf c}(k)$.
For matrices ${\bf A}$ and ${\bf B}$, with compatible dimensions,
${\bf A}\odot {\bf B}$ denotes the Khatri-Rao (columnwise Kronecker) product, ${\bf A}\circledast {\bf B}$ denotes the Hadamard
(element-wise) product, and ${\bf A}*{\bf B}$ denotes the matrix inner product, that is 
${\bf A}*{\bf B}:= {\rm trace}({\bf A}^T {\bf B}) = \sum_{i,j} {\bf A}_{i,j} {\bf B}_{i,j}$.

\section{non-negative tensor factorization}
\label{Section_NTF}

Let tensor $\underline{\bf X}^o\in\mathbb{R}^{I\times J\times K}_+$ admit a
non-negative\footnote{Note that, due to the non-negativity constraints on the latent factors,
$F$ can be higher than the rank of $\underline{\bf X}^o$.} CP decomposition of order $F$
\begin{equation*}
\underline{\bf X}^o  = [{\bf A}^o, {\bf B}^o, {\bf C}^o] =
\sum_{f=1}^F {\bf a}^o_f \circ {\bf c}^o_f \circ {\bf c}^o_f,
\label{Xo_def}
\end{equation*}
where ${\bf A}^o=[{\bf a}_1^o ~\cdots ~ {\bf a}^o_F]\in\mathbb{R}_+^{I\times F}$,
${\bf B}^o=[{\bf b}_1^o ~\cdots ~ {\bf b}^o_F]\in\mathbb{R}_+^{J\times F}$, 
and ${\bf C}^o=[{\bf c}_1^o ~\cdots ~ {\bf c}^o_F]\in\mathbb{R}_+^{K\times F}$.
We observe a noisy version of $\underline{\bf X}^o$ expressed as
\begin{equation*}
\underline{\bf X} = \underline{\bf X}^o + \underline{\bf E}.
\label{X_def}
\end{equation*}
In order to estimate ${\bf A}^o$, ${\bf B}^o$, and ${\bf C}^o$, we
compute matrices ${\bf A}\in\mathbb{R}_+^{I\times F}$, ${\bf B}\in\mathbb{R}_+^{J\times F}$, and
${\bf C}\in\mathbb{R}_+^{K\times F}$ that solve the optimization problem
\begin{equation}
\begin{array}{cl}
\displaystyle\min_{{\bf A}, {\bf B}, {\bf C}} & f_{\underline{\bf X}}({\bf A}, {\bf B}, {\bf C})  \cr
{\rm subject~to} & {\bf A}\ge  {\bf 0}, ~{\bf B}\ge {\bf 0}, ~{\bf C}\ge {\bf 0},
\end{array}
\label{Problem_NTF}
\end{equation}
where $f$ is a function measuring the quality of the factorization, ${\bf 0}$ is the zero matrix
of appropriate dimensions, and the inequalities are element-wise.
A common choice for $f_{\underline{\bf X}}$, motivated via maximum likelihood estimation for $\underline{\bf E}$ 
with Gaussian independent and identically distributed (i.i.d.) elements, is
\begin{equation}
f_{\underline{\bf X}}({\bf A}, {\bf B}, {\bf C}) =
\frac{1}{2} \, \left\| \underline{\bf X} - [{\bf A}, {\bf B}, {\bf C}] \right\|_F^2.
\end{equation}
Let $\underline{\bf W}=[{\bf A}, {\bf B}, {\bf C}]$ and ${\bf W}^{(1)}$, ${\bf W}^{(2)}$, and ${\bf W}^{(3)}$ be
the matrix unfoldings of $\underline{\bf W}$, with respect to the first, second, and third dimension,
respectively. Then,
\begin{equation}
\begin{split}
{\bf W}^{(1)} & = {\bf A} \, ({\bf C} \odot {\bf B})^T, \cr
{\bf W}^{(2)} & = {\bf B} \, ({\bf C} \odot {\bf A})^T, \cr
{\bf W}^{(3)} & = {\bf C} \, ({\bf B} \odot {\bf A})^T,
\end{split}
\end{equation}
and $f_{\underline{\bf X}}$ can be equivalently expressed as
\begin{equation}
\begin{split}
f_{\underline{\bf X}}({\bf A}, {\bf B}, {\bf C}) & = \frac{1}{2} \,\left\| {\bf X}^{(1)} - {\bf A} \, ({\bf C} \odot {\bf B})^T \right\|_F^2 \cr
& = \frac{1}{2} \,\left\| {\bf X}^{(2)}  - {\bf B} \, ({\bf C} \odot {\bf A})^T \right\|_F^2 \cr
& = \frac{1}{2} \,\left\|  {\bf X}^{(3)} - {\bf C} \, ({\bf B} \odot {\bf A})^T \right\|_F^2.
\end{split}
\label{f_X_matr}
\end{equation}
These expressions are the basis for ALS-type CP optimization, because they enable simple
linear least-squares updating of one matrix given the other two. Using NALS for each update step
is a popular approach for the solution of (\ref{Problem_NTF}), but non-negativity
brings a significant computational burden relative to plain ALS and also complicates the development of parallel
algorithms for NTF. It is worth noting that the above expressions will also prove useful during the development of the
ADMoM-based NTF algorithm.

\section{ADMoM}
\label{Section_ADMoM}

ADMoM is a technique for the solution of optimization problems of the form \cite{ADMM_Boyd_et_al_NOW}
\begin{equation}
\begin{array}{cl}
\displaystyle\min_{{\bf x}, {\bf z}} & f({\bf x}) + g({\bf z}) \cr
{\rm subject~to} & {\bf A} {\bf x} + {\bf B} {\bf z} = {\bf c},
\end{array}
\label{Problem_P_ADMM}
\end{equation}
where ${\bf x}\in\mathbb{R}^{n_1}$, ${\bf z}\in\mathbb{R}^{n_2}$, ${\bf A}\in\mathbb{R}^{m\times n_1}$,
${\bf B}\in\mathbb{R}^{m\times n_2}$, ${\bf c}\in\mathbb{R}^m$, $f:\mathbb{R}^{n_1}\rightarrow \mathbb{R}$,
and $g:\mathbb{R}^{n_2}\rightarrow \mathbb{R}$.

The augmented Lagrangian for problem (\ref{Problem_P_ADMM}) is
\begin{equation}
\begin{split}
L_{\rho}({\bf x}, {\bf z}, {\bf y}) & = f({\bf x}) + g({\bf z}) + {\bf y}^T ({\bf A} {\bf x} + {\bf B} {\bf z}-{\bf c}) \cr
& \hspace{1.5cm} + \frac{\rho}{2} \|{\bf A}{\bf x}+{\bf B}{\bf z}-{\bf c}\|_2^2,
\label{L_rho_Problem_P_ADMM}
\end{split}
\end{equation}
where $\rho>0$ is a penalty parameter.
Assuming that at time instant $k$ we have computed ${\bf z}^k$ and ${\bf y}^k$, which comprise the state of
the algorithm, the $(k+1)$-st iteration of ADMoM is\footnote{ Note that ${\bf y}^{kT}$ is 
shorthand notation for $\left({\bf y}^k\right)^T$.}
\begin{equation*}
\begin{split}
{\bf x}^{k+1} & = \underset{\bf x} {\rm argmin}\left( f({\bf x}) + {\bf y}^{kT} {\bf A}{\bf x} +
\frac{\rho}{2} \|{\bf A} {\bf x} +{\bf B} {\bf z}^k -{\bf c} \|_2^2 \right) \cr
{\bf z}^{k+1} & = \underset{\bf z} {\rm argmin}\left(g({\bf z}) + {\bf y}^{kT} {\bf B}{\bf z} +
\frac{\rho}{2} \| {\bf A} {\bf x}^{k+1} +{\bf B} {\bf z} -{\bf c} \|_2^2    \right) \cr
{\bf y}^{k+1} & = {\bf y}^k + \rho \, ({\bf A}{\bf x}^{k+1} +{\bf B} {\bf z}^{k+1}-{\bf c}).
\end{split}
\end{equation*}
It can be shown that, under certain conditions (among them convexity of $f$ and $g$), ADMoM
converges in a certain sense (see \cite{ADMM_Boyd_et_al_NOW} for an excellent review of
ADMoM, including some convergence analysis results). However, ADMoM can be used even when problem
(\ref{Problem_P_ADMM}) is {\em non-convex}. In this case, we use ADMoM with the goal of reaching
a good local minimum \cite{ADMM_Boyd_et_al_NOW}. Note that this is all we can realistically hope for anyway,
irrespective of approach or algorithm used, since tensor factorization is NP-hard \cite{tensorNP}.

\subsection{ADMoM for set-constrained optimization}
\label{subsection_ADMoM_set_constraint}

Let us consider the set-constrained optimization problem
\begin{equation*}
\begin{array}{cl}
\displaystyle\min_{{\bf x}} & f({\bf x}) \cr
{\rm subject~to} & {\bf x}\in{\cal X},
\end{array}
\label{Problem_P_X}
\end{equation*}
where ${\bf x}\in\mathbb{R}^n$ and ${\cal X}\subseteq\mathbb{R}^n$ is a closed convex set.
At first sight, this problem does not seem suitable for ADMoM. However,
if we introduce variable ${\bf z} \in\mathbb{R}^n$,
we can write the equivalent problem
\begin{equation}
\begin{array}{cl}
\displaystyle\min_{{\bf x}, {\bf z}} & f({\bf x}) + g({\bf z}) \cr
{\rm subject~to} & {\bf x}-{\bf z}=0,
\end{array}
\label{Problem_P_ADMM_form}
\end{equation}
where $g$ is the indicator function of set ${\cal X}$, that is,
\begin{equation*}
g({\bf z}) :=
\left\{ \begin{array}{ll}  0, & {\bf z}\in{\cal X}, \cr \infty, & {\bf z}\notin{\cal X}.\end{array}  \right.
\end{equation*}
Then it becomes clear that (\ref{Problem_P_ADMM_form}) can be solved via ADMoM.
Assuming that at time instant $k$ we have computed ${\bf z}^k$ and ${\bf y}^k$,
the $(k+1)$-st iteration of ADMoM is \cite{ADMM_Boyd_et_al_NOW}
\begin{equation*}
\begin{split}
{\bf x}^{k+1} & = \underset{\bf x} {\rm argmin}\left( f({\bf x}) + {\bf y}^{kT} {\bf x} +
\frac{\rho}{2} \| {\bf x} - {\bf z}^k \|_2^2 \right) \cr
{\bf z}^{k+1} & = \underset{\bf z} {\rm argmin}\left(g({\bf z}) - {\bf y}^{kT} {\bf z} +
\frac{\rho}{2} \|{\bf x}^{k+1} - {\bf z} \|_2^2    \right) \cr
& = \Pi_{\cal X}\left({\bf x}^{k+1}+ \frac{1}{\rho}{\bf y}^k \right) \cr
{\bf y}^{k+1} & = {\bf y}^k + \rho \, ({\bf x}^{k+1} - {\bf z}^{k+1}),
\end{split}
\end{equation*}
where $\Pi_{\cal X}$ denotes projection (in the Euclidean norm) onto ${\cal X}$.

\section{ADMoM for NTF}
\label{Section_ADMoM_NTF}

In this section, we adopt the approach of subsection \ref{subsection_ADMoM_set_constraint}
and develop an ADMoM-based NTF algorithm. At first, we must put the NTF problem (\ref{Problem_NTF}) into ADMoM form.
Towards this end, we introduce auxiliary variables $\tilde{\bf A}\in\mathbb{R}^{I\times F}$,
$\tilde{\bf B}\in\mathbb{R}^{J\times F}$, and $\tilde{\bf C}\in\mathbb{R}^{K \times F}$ and consider the
equivalent optimization problem
\begin{equation}
\begin{array}{cl}
\displaystyle\min_{{\bf A}, \tilde{\bf A}, {\bf B},\tilde{\bf B}, {\bf C}, \tilde{\bf C}} &
f_{\underline{\bf X}}({\bf A}, {\bf B}, {\bf C}) +g(\tilde{\bf A}) + g(\tilde{\bf B}) + g(\tilde{\bf C}) \cr
{\rm subject~to} &  \hspace{-.2cm}{\bf A}-\tilde{\bf A} =
{\bf 0}, \,{\bf B}-\tilde{\bf B} = {\bf 0}, \,{\bf C}-\tilde{\bf C} = {\bf 0},
\end{array}
\label{Problem_NTF_ADMM_1}
\end{equation}
where, for any matrix argument ${\bf M}$,
\begin{equation}
g({\bf M}) := \left\{ \begin{array}{ll} 0, & {\rm if}~{\bf M} \ge {\bf 0}, \cr
\infty, & {\rm otherwise}.  \end{array} \right.
\label{g_M_def}
\end{equation}
We introduce the dual variables ${\bf Y}_{\bf A}\in\mathbb{R}^{I\times F}$,
${\bf Y}_{\bf B}\in\mathbb{R}^{J \times F}$, and ${\bf Y}_{\bf C}\in\mathbb{R}^{K \times F}$, and the vector of
penalty terms $\mbox{\boldmath{$\rho$}}:=[\rho_{\bf A}~\rho_{\bf B}~\rho_{\bf C}]^T$.
The augmented Lagrangian is given in (\ref{L_rho_NTF}), at the top of this page.

\begin{table*}
\normalsize
\begin{equation}
\begin{split}
L_{\mbox{\small\boldmath{$\rho$}}}({\bf A}, {\bf B}, {\bf C}, \tilde{\bf A}, \tilde{\bf B}, \tilde{\bf C},
{\bf Y}_{\bf A}, {\bf Y}_{\bf B}, {\bf Y}_{\bf C}) & =
f_{\underline{\bf X}}({\bf A}, {\bf B}, {\bf C}) + g(\tilde{\bf A}) + g(\tilde{\bf B}) +  g(\tilde{\bf C}) \cr
& + {\bf Y}_{\bf A} * ({\bf A}-\tilde{\bf A})  + \frac{\rho_{\bf A}}{2} \, \| {\bf A}-\tilde{\bf A} \|_F^2  \cr
& + {\bf Y}_{\bf B} * ({\bf B}-\tilde{\bf B})  + \frac{\rho_{\bf B}}{2} \, \|{\bf B}-\tilde{\bf B}\|_F^2\cr
& + {\bf Y}_{\bf C} * ({\bf C}-\tilde{\bf C})  + \frac{\rho_{\bf C}}{2} \, \|{\bf C}-\tilde{\bf C}\|_F^2.
\end{split}
\label{L_rho_NTF}
\end{equation}
\hrulefill
\end{table*}

The ADMoM for problem (\ref{Problem_NTF_ADMM_1}) is as follows:
\begin{equation}
\begin{split}
& ({\bf A}^{k+1}, {\bf B}^{k+1}, {\bf C}^{k+1})  = \underset{{\bf A}, {\bf B}, {\bf C}}
{\rm argmin} \left( f_{\underline{\bf X}}({\bf A},{\bf B}, {\bf C}) \right. \cr
& \qquad \qquad \qquad \qquad + {\bf Y}^k_{\bf A} * {\bf A} + \frac{\rho_{\bf A}}{2} \,\|{\bf A}-\tilde{\bf A}^k \|_F^2   \cr
& \qquad \qquad \qquad \qquad + {\bf Y}^k_{\bf B} * {\bf B} + \frac{\rho_{\bf B}}{2} \,\|{\bf B}-\tilde{\bf B}^k\|_F^2 \cr
& \qquad \qquad \qquad \qquad \left. + {\bf Y}^k_{\bf C} * {\bf C} + \frac{\rho_{\bf C}}{2} \, \|{\bf C}-\tilde{\bf C}^k\|_F^2  \right)\cr
& \qquad \tilde{\bf A}^{k+1}  = \left( {\bf A}^{k+1} + \frac{1}{\rho_{\bf A}} {\bf Y}_{\bf A}^k \right)_+ \cr
& \qquad \tilde{\bf B}^{k+1}  = \left( {\bf B}^{k+1} + \frac{1}{\rho_{\bf B}} {\bf Y}_{\bf B}^k \right)_+ \cr
& \qquad \tilde{\bf C}^{k+1}  = \left( {\bf C}^{k+1} + \frac{1}{\rho_{\bf C}} {\bf Y}_{\bf C}^k \right)_+ \cr
& \qquad {\bf Y}_{\bf A}^{k+1}  = {\bf Y}_{\bf A}^{k} + \rho_{\bf A} \left( {\bf A}^{k+1}- \tilde{\bf A}^{k+1} \right) \cr
& \qquad {\bf Y}_{\bf B}^{k+1}  = {\bf Y}_{\bf B}^{k} + \rho_{\bf B} \left( {\bf B}^{k+1}- \tilde{\bf B}^{k+1} \right) \cr
& \qquad {\bf Y}_{\bf C}^{k+1}  = {\bf Y}_{\bf C}^{k} + \rho_{\bf C} \left( {\bf C}^{k+1}- \tilde{\bf C}^{k+1} \right).
\end{split}
\label{ADMoM_NTF}
\end{equation}

\begin{table*}
\normalsize
\begin{equation}
\begin{split}
{\bf A}^{k+1} & = \underset{\bf A} {\rm argmin} \left( \frac{1}{2} \, \|{\bf X}^{(1)} - {\bf A}  ({\bf C}^k \odot {\bf B}^k)^T \|_F^2
+ {\bf Y}^k_{\bf A} * {\bf A} + \frac{\rho_{\bf A}}{2} \, \|{\bf A}-\tilde{\bf A}^k\|_F^2 \right) \cr
& = \left( {\bf X}^{(1)} ({\bf C}^k \odot {\bf B}^k) + \rho_{\bf A} \tilde{\bf A}^{k} - {\bf Y}_{\bf A}^{k} ) \right)
\left( ({\bf C}^{k} \odot {\bf B}^k)^T ({\bf C}^{k} \odot {\bf B}^k)   + \rho_{\bf A} {\bf I}_{F} \right)^{-1}, \cr
{\bf B}^{k+1} & = \underset{\bf B} {\rm argmin} \left( \frac{1}{2} \, \|{\bf X}^{(2)} - {\bf B}  ({\bf C}^k \odot {\bf A}^{k+1})^T \|_F^2 +
{\bf Y}^k_{\bf B} * {\bf B} + \frac{\rho_{\bf B}}{2} \, \|{\bf B}-\tilde{\bf B}^k \|_F^2 \right) \cr
& = \left( {\bf X}^{(2)} ({\bf C}^k \odot {\bf A}^{k+1}) + \rho_{\bf B} \tilde{\bf B}^{k} - {\bf Y}_{\bf B}^{k} \right)
\left( ({\bf C}^{k} \odot {\bf A}^{k+1})^T ({\bf C}^{k} \odot {\bf A}^{k+1})   + \rho_{\bf B} {\bf I}_{F} \right)^{-1}, \cr
{\bf C}^{k+1} & = \underset{\bf C} {\rm argmin} \left( \frac{1}{2} \, \|{\bf X}^{(3)} - {\bf C}  ({\bf B}^{k+1} \odot {\bf A}^{k+1})^T \|_F^2 +
{\bf Y}^k_{\bf C} * {\bf C} + \frac{\rho_{\bf C}}{2} \, \|{\bf C}-\tilde{\bf C}^k\|_F^2 \right) \cr
& = \left( {\bf X}^{(3)} ({\bf B}^{k+1} \odot {\bf A}^{k+1}) + \rho_{\bf C} \tilde{\bf C}^{k} - {\bf Y}_{\bf C}^{k} \right)
\left( ({\bf B}^{k+1} \odot {\bf A}^{k+1})^T ({\bf B}^{k+1} \odot {\bf A}^{k+1})   + \rho_{\bf C} {\bf I}_{F} \right)^{-1}.
\end{split}
\label{ADMM_NTF_ALS_C}
\end{equation}
\hrulefill
\end{table*}

\noindent
The minimization problem in the first line of (\ref{ADMoM_NTF}) is non-convex.
Using the equivalent expressions for $f_{\underline{\bf X}}$ in (\ref{f_X_matr}), we propose the
alternating optimization scheme of (\ref{ADMM_NTF_ALS_C}), at the top of the next page.
The updates of (\ref{ADMM_NTF_ALS_C}) can be executed either for a predetermined number of iterations, or until
convergence.\footnote{In our implementations, we execute these updates once per ADMoM iteration.}
We observe that, during each ADMoM iteration, we avoid the solution of constrained optimization problems.
This seems favorable, especially in the cases where the size of the problem is (very) large.

We note that ${\bf A}^k$, ${\bf B}^k$, and ${\bf C}^k$ are {\em not} necessarily non-negative.
They become non-negative (or, at least, their negative elements become very small) upon convergence.
On the other hand, $\tilde{\bf A}^k$, $\tilde{\bf B}^k$, and $\tilde{\bf C}^k$ are by construction
non-negative.

\subsection{Computational complexity per iteration}

Each ADMoM iteration consists of simple matrix operations. Thus, rough estimates of its computational complexity can
be easily derived (of course, accurate estimates can be derived after fixing the algorithms
that implement the matrix operations).

A rough estimate for the computational complexity of the update of
${\bf A}^k$ (see the first update in (\ref{ADMM_NTF_ALS_C})) can be derived as follows:
\begin{enumerate}
\item $O(n F+IF)$ for the computation of the term
${\bf X}^{(1)} ({\bf C}^k \odot {\bf B}^k) + \rho_{\bf A} \tilde{\bf A}^{k} - {\bf Y}_{\bf A}^{k}$, 
where $n$ is the number of nonzero elements of tensor $\underline{\bf X}$ (also of matrix ${\bf X}^{(1)}$). 
Note that $n=IJK$ for dense tensors, but for sparse tensors $n \ll IJK$.  This is because the 
product ${\bf X}^{(1)} ({\bf C}^k \odot {\bf B}^k)$ can be computed with $3nF$ flops, by exploiting 
sparsity and the structure of the Khatri-Rao product \cite{BaderKolda2007,kolda2008scalable,TensorToolbox}. 
With $5nF$ flops, it is possible to parallelize this computation \cite{kang2012gigatensor}. Efficient 
(in terms of favorable memory access pattern) in-place computation of all three products needed for the 
update of ${\bf A}^k$, ${\bf B}^k$, ${\bf C}^k$ from a single copy of $\underline{\bf X}$ has been recently 
considered in \cite{RavSidSmiKar:Asilo2014}, which also features potential flop gains as a side-benefit.
\item $O((K+J)F^2)$ for the computation of the term
$({\bf C}^{k} \odot {\bf B}^k)^T ({\bf C}^{k} \odot {\bf B}^k)   + \rho_{\bf A} {\bf I}_{F}$, and $O(F^3)$ 
for its Cholesky decomposition. This is because 
$({\bf C}^{k} \odot {\bf B}^k)^T ({\bf C}^{k} \odot 
{\bf B}^k)$ $=$ $\left(\left({\bf C}^{k}\right)^T {\bf C}^{k}\right) \circledast \left(\left({\bf B}^{k}\right)^T {\bf B}^{k}\right)$.
\item $O(F^2I)$ for the computation of the system solution that gives the updated value ${\bf A}^{k+1}$.
\end{enumerate}
Analogous estimates can be derived for the updates of ${\bf B}^k$ and ${\bf C}^k$. Finally,
the updates of the auxiliary and dual variables require, in total,
$O\left((I+J+K)F\right)$ arithmetic operations.

\subsection{Convergence}
\label{subsection_convergence}

Let $\boldsymbol{Z}:=({\bf A}, {\bf B}, {\bf C}, \tilde{\bf A}, \tilde{\bf B}, \tilde{\bf C}, {\bf Y}_{\bf A},
{\bf Y}_{\bf B},  {\bf Y}_{\bf C})$. It can be proven that $\boldsymbol{Z}$ is a Karush-Kuhn-Tucker (KKT)
point for the NTF problem (\ref{Problem_NTF_ADMM_1}) if
\begin{equation}
\begin{split}
& \left({\bf X}^{(1)} - {\bf A} ({\bf C}\odot {\bf B})^T \right) ({\bf C}\odot {\bf B}) - {\bf Y}_{\bf A}  = {\bf 0} \cr
& \left({\bf X}^{(2)} - {\bf B} ({\bf C}\odot {\bf A})^T \right) ({\bf C}\odot {\bf A}) - {\bf Y}_{\bf B}  = {\bf 0} \cr
& \left({\bf X}^{(3)} - {\bf C} ({\bf B}\odot {\bf A})^T \right) ({\bf B}\odot {\bf A}) - {\bf Y}_{\bf C}  = {\bf 0} \cr
& {\bf A} - \tilde{\bf A} = {\bf 0}, ~ {\bf B} - \tilde{\bf B} = {\bf 0}, ~ {\bf C} - \tilde{\bf C} = {\bf 0} \cr
& {\bf Y}_{\bf A}\le {\bf 0},~ {\bf Y}_{\bf B}\le {\bf 0},~ {\bf Y}_{\bf C}  \le {\bf 0} \cr
& {\bf Y}_{\bf A} \circledast \tilde{\bf A}  = {\bf 0}, ~ {\bf Y}_{\bf B} \circledast \tilde{\bf B}  = {\bf 0}, ~
{\bf Y}_{\bf C} \circledast \tilde{\bf C} = {\bf 0}.
\end{split}
\label{KKT_NTF}
\end{equation}

\begin{proposition}
Let $\{\boldsymbol{Z}^k\}$ be a sequence generated by ADMoM for NTF that satisfies condition
\begin{equation}
\lim_{k\rightarrow \infty} (\boldsymbol{Z}^{k+1} - \boldsymbol{Z}^k) = {\bf 0}.
\end{equation}
Then, any accumulation point of $\{\boldsymbol{Z}^k\}$ is a KKT point of problem (\ref{Problem_NTF_ADMM_1}).
Consequently, any accumulation point of $\{{\bf A}^k, {\bf B}^k, {\bf C}^k\}$ is a KKT point of problem (\ref{Problem_NTF}).
\label{proposition_1}
\end{proposition}

{\em Proof:} The proof follows closely the steps of the proof of Proposition 2.1 of
\cite{Xu_et_al_2010} and is omitted.\footnote{See report \cite{Liavas_Sidiropoulos_Report} for a detailed proof.} \hfill$\Box$

Proposition \ref{proposition_1} implies that, whenever $\{\boldsymbol{Z}^k\}$ converges, it converges to a KKT point.
We will further discuss ADMoM convergence from a practical point of view in the section with the numerical experiments.

\subsection{Stopping criteria}
\label{subsection_stopping_criterion}

The primal residual for variable ${\bf A}^k$ is defined as
\begin{equation}
{\bf P}_{\bf A}^k := {\bf A}^k - \tilde{\bf A}^k,
\label{R_Mk_def}
\end{equation}
while quantity
\begin{equation}
{\bf D}_{\tilde{\bf A}}^k := \rho_{\bf A} (\tilde{\bf A}^k - \tilde{\bf A}^{k-1})
\label{D_Mk_def}
\end{equation}
can be viewed as a dual feasibility residual (see \cite[Section 3.3]{ADMM_Boyd_et_al_NOW}).
We analogously define ${\bf P}_{\bf B}^k$, ${\bf D}_{\tilde{\bf B}}^k$, ${\bf P}_{\bf C}^k$, and
${\bf D}_{\tilde{\bf C}}^k$.

We stop the algorithm if all primal and dual residuals are sufficiently small. More specifically,
we introduce small positive constants $\epsilon^{\rm abs}$ and $\epsilon^{\rm rel}$
and consider ${\bf P}_{\bf A}^k$ and ${\bf D}_{\tilde{\bf A}}^k$ small if
\begin{equation}
\| {\bf P}_{\bf A}^k\|_F \le \sqrt{IF}\,\epsilon^{\rm abs} + \epsilon^{\rm rel} \max\left\{ \|{\bf A}^k\|_F, \|\tilde{\bf A}^k\|_F \right\},
\label{stopping_cond_prim_A_def}
\end{equation}
\begin{equation}
\|{\bf D}_{\tilde{\bf A}}^k\|_F \le \sqrt{IF}\, \epsilon^{\rm abs} + \epsilon^{\rm rel} \, \|{\bf Y}_{\bf A}^k\|_F.
\label{stopping_cond_dual_A_def}
\end{equation}
Analogous conditions apply for the other residuals. Reasonable values for $\epsilon^{\rm rel}$ are $\epsilon^{\rm rel}\lessapprox 10^{-3}$,
while the value of $\epsilon^{\rm abs}$ depends on the scale of the values of the latent factors.

We note that stopping criteria (\ref{stopping_cond_prim_A_def}) and (\ref{stopping_cond_dual_A_def})
involve quantities of the size of the latent factors which, in most cases, is small compared to the size of the tensor.
Thus, their computation, even during every ADMoM iteration, is not computationally demanding.

\subsection{Varying penalty parameters}
\label{subsection_varying_penalty}

We have found very useful in practice to vary the values of each one of the 
penalty parameters, $\rho_{\bf A}$, $\rho_{\bf B}$, and
$\rho_{\bf C}$, depending on the size of the corresponding primal and dual residuals 
(see \cite[Section 3.4]{ADMM_Boyd_et_al_NOW}). More specifically,
the penalty parameters $\rho_{\bf M}^k$, for ${\bf M}={\bf A}, {\bf B}, {\bf C}$, are updated as follows:
\begin{equation}
\rho_{\bf M}^{k+1} = \left\{
\begin{array}{ll}
\tau^{\rm incr} \rho_{\bf M}^{k}, & \mbox{if}~ \|{\bf P}_{\bf M}^k\|_F > \mu \,\|{\bf D}_{\tilde{\bf M}}^k\|_F, \cr
\rho_{\bf M}^{k}/\tau^{\rm decr},  & \mbox{if}~ \|{\bf D}_{\tilde{\bf M}}^k\|_F > \mu\, \|{\bf P}_{\bf M}^k\|_F, \cr
\rho_{\bf M}^{k}, & \mbox{otherwise},
\end{array}
\right.
\label{varying_penalty_rule}
\end{equation}
where $\mu>1$, $\tau^{\rm incr}>1$, and $\tau^{\rm decr}>1$ are the adaptation parameters.
Large values of $\rho_{\bf M}$ place large penalty on violations of primal feasibility, 
leading to small primal residuals, while small values of $\rho_{\bf M}$ tend to reduce the dual residuals.

\subsection{ADMoM for tensor factorization with structural constraints}
\label{subsection_structural_constraints}

ADMoM can easily handle certain structural constraints on the latent factors
\cite{ADMM_Boyd_et_al_NOW}, \cite{Xu_et_al_2013}. For example, if we want to solve an NTF problem with the
added constraint that the number of nonzero elements of ${\bf A}$ is lower than or equal to a given
number $c_{\bf A}$, then we can adopt an approach similar to that followed in Section \ref{Section_ADMoM_NTF}
with the only difference being that, instead of using $g(\tilde{\bf A})$ defined in
(\ref{g_M_def}), we use $g_{c_{\bf A}}(\tilde{\bf A})$ where, for any matrix argument ${\bf M}$,
\begin{equation}
g_c({\bf M}) := \left\{ \begin{array}{ll} 0, & \mbox{if}~ {\bf M}\ge {\bf 0}~\mbox{and}~\|{\bf M}\|_0 \le c, \cr
\infty, & \mbox{otherwise}. \end{array} \right.
\label{g_new_def}
\end{equation}
The only difference between the ADMoM for this case and the one presented in (\ref{ADMoM_NTF}) and (\ref{ADMM_NTF_ALS_C})
is in the update of $\tilde{\bf A}^k$. More specifically, instead of using projection
onto the set of non-negative matrices, we must use projection onto the set of non-negative matrices with at most
$c_{\bf A}$ nonzero elements, which can be easily computed through sorting of the
elements of $\tilde{\bf A}^k$.
Using analogous arguments, we can incorporate into our ADMoM framework box or other set constraints on the latent factors.
The development of the corresponding ADMoM is almost trivial if  projection onto the constraint
set is easy.  

Thorough study of ADMoM-based algorithms for tensor factorization and/or completion with more complicated structural
constraints is a topic of future research.

\section{Distributed ADMoM for large NTF}
\label{Section_ADMoM_LNTF}

In this section, we assume that all dimensions of tensor $\underline{\bf X}$ are large and derive
an ADMoM-based NTF that is suitable for parallel implementation. Of course, our framework can handle
the cases where only one or two of the dimensions of
$\underline{\bf X}$ are large. 

\subsection{Matrix unfoldings in terms of partitioned matrix factors}
\label{subsection_partitioned_matrix_unfoldings}

Let $\underline{\bf W}=[{\bf A}, {\bf B}, {\bf C}]$, and ${\bf A}$, ${\bf B}$, and ${\bf C}$ be partitioned as
\begin{equation}
{\bf A}= \left[\begin{array}{c}{\bf A}_1 \cr \vdots \cr {\bf A}_{N_A}\end{array}\right], ~
{\bf B}= \left[\begin{array}{c}{\bf B}_1 \cr \vdots \cr {\bf B}_{N_B}\end{array}\right], ~
{\bf C}= \left[\begin{array}{c}{\bf C}_1 \cr \vdots \cr {\bf C}_{N_C}\end{array}\right],
\end{equation}
with ${\bf A}_{n_A}\in\mathbb{R}^{I_{n_A}\times F}$, for $n_A=1,\ldots,N_A$,
$\sum_{n_A=1}^{N_A}I_{n_A}=I$,
${\bf B}_{n_B}\in\mathbb{R}^{J_{n_B}\times F}$, for $n_B=1,\ldots,N_B$, $\sum_{n_B=1}^{N_B}J_{n_B}=J$, and
${\bf C}_{n_C}\in\mathbb{R}^{K_{n_C}\times F}$, for $n_C=1,\ldots,N_C$, $\sum_{n_C=1}^{N_C} K_{n_C}=K$.

We first derive partitionings of the matrix unfoldings of
$\underline{\bf W}$ in terms of (the blocks of) matrices ${\bf A}$, ${\bf B}$, and ${\bf C}$.
Towards this end, we write ${\bf W}^{(1)}$ as
\begin{equation*}
\begin{split}
{\bf W}^{(1)} & = {\bf A} ({\bf C}  \odot {\bf B})^T \cr
& = \left[  \begin{array}{c} {\bf A}_1 \cr \vdots \cr {\bf A}_{N_A} \end{array} \right]
\left(  \left[  \begin{array}{c} {\bf C}_1 \cr \vdots \cr {\bf C}_{N_C}  \end{array} \right]  \odot
{\bf B}  \right)^T \cr
& =  \left[  \begin{array}{c} {\bf A}_1 \cr \vdots \cr {\bf A}_{N_A} \end{array} \right]
\left(  \left[  \begin{array}{c} {\bf C}_1 \odot {\bf B} \cr \vdots \cr {\bf C}_{N_C}\odot {\bf B}  \end{array} \right]  \right)^T \cr
& = \left[  \begin{array}{c} {\bf A}_1 \cr \vdots \cr {\bf A}_{N_A} \end{array} \right]
\left(  \left[ \begin{array}{ccc} ({\bf C}_1 \odot {\bf B})^T & \cdots & ({\bf C}_{N_C}\odot {\bf B})^T  \end{array} \right]  \right).
\end{split}
\end{equation*}
Thus, ${\bf W}^{(1)}$ can be partitioned as
\begin{equation*}
\begin{split}
{\bf W}^{(1)} & =
\left[ \begin{array}{ccc}
{\bf W}^{(1)}_{1,1} & \cdots & {\bf W}^{(1)}_{1,N_C} \cr
\vdots & \ddots & \vdots \cr
{\bf W}^{(1)}_{N_A,1} & \cdots & {\bf W}^{(1)}_{N_A,N_C}
\end{array}\right],
\end{split}
\end{equation*}
where the $(n_A,n_C)$-th block of ${\bf W}^{(1)}$  is equal to the
$I_{n_A} \times (J K_{n_C})$ matrix ${\bf W}^{(1)}_{n_A,n_C}={\bf A}_{n_A} ({\bf C}_{n_C} \odot {\bf B})^T$,
for $n_A=1,\ldots,N_A$ and $n_C=1,\ldots,N_C$.

Similarly, it can be shown that ${\bf W}^{(2)}$ can be partitioned into blocks
${\bf W}^{(2)}_{n_B, n_C}={\bf B}_{n_B} ({\bf C}_{n_C} \odot {\bf A})^T$,
of dimensions $J_{n_B} \times (I K_{n_C})$,
for $n_B=1,\ldots,N_B$ and $n_C=1,\ldots,N_C$,
and ${\bf W}^{(3)}$ can be partitioned into blocks
${\bf W}^{(3)}_{n_C, n_B} = {\bf C}_{n_C}  ({\bf B}_{n_B} \odot {\bf A})^T$,
of dimensions $K_{n_C} \times (I J_{n_B})$, for $n_C=1,\ldots,N_C$ and
$n_B=1,\ldots,N_B$.\footnote{An extension of the above partitioning scheme to higher
order tensors appears in Appendix \ref{partitionings_extension}.}

If we partition ${\bf X}^{(1)}$, ${\bf X}^{(2)}$, and ${\bf X}^{(3)}$ accordingly, then we can write
\begin{equation}
\begin{split}
& f_{\underline{\bf X}}({\bf A}, {\bf B}, {\bf C})  =  \sum_{n_A=1}^{N_A} \sum_{n_C=1}^{N_C} \frac{1}{2} \,
\|{\bf X}^{(1)}_{n_A,n_C} - {\bf A}_{n_A} ( {\bf C}_{n_C} \odot {\bf B})^T \|_F^2  \cr
& \quad =  \sum_{n_B=1}^{N_B} \sum_{n_C=1}^{N_C} \frac{1}{2} \,
\|{\bf X}^{(2)}_{n_B,n_C} - {\bf B}_{n_B} ( {\bf C}_{n_C} \odot {\bf A})^T \|_F^2  \cr
& \quad =  \sum_{n_C=1}^{N_C} \sum_{n_B=1}^{N_B} \frac{1}{2} \,
\|{\bf X}^{(3)}_{n_C,n_B} - {\bf C}_{n_C} ( {\bf B}_{n_B} \odot {\bf A})^T \|_F^2.
\end{split}
\label{f_X_vl}
\end{equation}
These expressions will be fundamental for the development of the distributed ADMoM for large NTF.

\subsection{Distributed ADMoM for large NTF}
\label{subsection_ADMOM_vl}

In order to put the large NTF problem into ADMoM form, we introduce auxiliary variables
$\tilde{\bf A}=[\tilde{\bf A}_1^T~\cdots ~\tilde{\bf A}_{N_A}^T]^T$, with
$\tilde{\bf A}_{n_A}\in\mathbb{R}^{I_{n_A}\times F}$, for $n_A=1,\ldots,N_A$,
$\tilde{\bf B}=[\tilde{\bf B}_1^T~\cdots ~\tilde{\bf B}_{N_B}^T]^T$, with
$\tilde{\bf B}_{n_B}\in\mathbb{R}^{J_{n_B}\times F}$, for $n_B=1,\ldots,N_B$,
and
$\tilde{\bf C}=[\tilde{\bf C}_1^T~\cdots ~\tilde{\bf C}_{N_C}^T]^T$, with
$\tilde{\bf C}_{n_C}\in\mathbb{R}^{K_{n_C}\times F}$, for $n_C=1,\ldots,N_C$,
and consider the equivalent problem
\begin{equation}
\begin{array}{cl}
\displaystyle\min_{{\bf A}, \tilde{\bf A}, {\bf B},\tilde{\bf B}, {\bf C}, \tilde{\bf C}} &
f_{\underline{\bf X}}({\bf A}, {\bf B}, {\bf C})
+ \sum_{n_A=1}^{N_A} g(\tilde{\bf A}_{n_A}) \cr
& \quad + \sum_{n_B=1}^{N_B} g(\tilde{\bf B}_{n_B}) + \sum_{n_C=1}^{N_C} g(\tilde{\bf C}_{n_C}) \cr
{\rm subject~to} & {\bf A}_{n_A}-\tilde{\bf A}_{n_A}={\bf 0}, \quad n_A=1,\ldots,N_A, \cr
& {\bf B}_{n_B} - \tilde{\bf B}_{n_B}={\bf 0}, \quad n_B=1,\ldots,N_B, \cr
& {\bf C}_{n_C}-\tilde{\bf C}_{n_C}={\bf 0}, \quad n_C=1,\ldots,N_C.
\end{array}
\label{Problem_partitioned_NTF_ADMM_vl}
\end{equation}
If we introduce dual variables ${\bf Y}_{\bf A} = [{\bf Y}_{{\bf A}_1}^T ~\cdots ~{\bf Y}_{{\bf A}_{N_A}}^T]^T$,
with ${\bf Y}_{{\bf A}_{n_A}}\in\mathbb{R}^{I_{n_A}\times F}$,
for $n_A=1,\ldots,N_A$,
${\bf Y}_{\bf B} = [{\bf Y}_{{\bf B}_1}^T ~\cdots ~{\bf Y}_{{\bf B}_{N_B}}^T]^T$, with
${\bf Y}_{{\bf B}_{n_B}}\in\mathbb{R}^{J_{n_B} \times F}$, for $n_B=1,\ldots,N_B$,  and
${\bf Y}_{\bf C} = [{\bf Y}_{{\bf C}_1}^T ~\cdots ~{\bf Y}_{{\bf C}_{N_A}}^T]^T$,
with ${\bf Y}_{{\bf C}_{n_C}} \in\mathbb{R}^{K_{n_C} \times F}$, for $n_C=1,\ldots,N_C$,
the augmented Lagrangian is written as in (\ref{L_rho_vl}), at the top of the next page.

\begin{table*}
\normalsize
\begin{equation}
\begin{split}
L_{\mbox{\small\boldmath{$\rho$}}}({\bf A}, {\bf B}, {\bf C}, \tilde{\bf A}, \tilde{\bf B},
\tilde{\bf C}, {\bf Y}_{\bf A}, {\bf Y}_{\bf B}, {\bf Y}_{{\bf C}})
& = f_{\underline{\bf X}}({\bf A}, {\bf B}, {\bf C}) + \sum_{n_A=1}^{N_A} g(\tilde{\bf A}_{n_A}) +
\sum_{n_B=1}^{N_B} g(\tilde{\bf B}_{n_B}) + \sum_{n_C=1}^{N_C} g(\tilde{\bf C}_{n_C})
\cr
& +  \sum_{n_A=1}^{N_A} \left( {\bf Y}_{{\bf A}_{n_A}} * ({\bf A}_{n_A} - \tilde{\bf A}_{n_A}) + \frac{\rho_{\bf A}}{2}\, \|{\bf A}_{n_A}-\tilde{\bf A}_{n_A}\|_F^2 \right) \cr
& +  \sum_{n_B=1}^{N_B} \left( {\bf Y}_{{\bf B}_{n_B}} * ({\bf B}_{n_B} - \tilde{\bf B}_{n_B}) + \frac{\rho_{\bf B}}{2}\, \|{\bf B}_{n_B}-\tilde{\bf B}_{n_B}\|_F^2 \right)\cr
& +  \sum_{n_C=1}^{N_C} \left( {\bf Y}_{{\bf C}_{n_C}} * ({\bf C}_{n_C} - \tilde{\bf C}_{n_C}) + \frac{\rho_{\bf C}}{2}\, \|{\bf C}_{n_C}-\tilde{\bf C}_{n_C} \|_F^2 \right).
\end{split}
\label{L_rho_vl}
\end{equation}
\hrulefill
\end{table*}

The ADMoM for this problem is as follows:
\begin{equation}
\begin{split}
& ({\bf A}^{k+1}, {\bf B}^{k+1}, {\bf C}^{k+1})  = \underset{{\bf A}, {\bf B}, {\bf C}}
{\rm argmin} \bigg( f_{\underline{\bf X}}({\bf A}, {\bf B}, {\bf C}) \Bigr. \cr
& \quad \qquad +   \sum_{n_A=1}^{N_A} \left( {\bf Y}^k_{{\bf A}_{n_A}} * {\bf A}_{n_A} + \frac{\rho_{\bf A}}{2} \, \|{\bf A}_{n_A}-\tilde{\bf A}^k_{n_A} \|_F^2 \right)  \cr
&  \quad \qquad +  \sum_{n_B=1}^{N_B} \left( {\bf Y}^k_{{\bf B}_{n_B}} * {\bf B}_{n_B} + \frac{\rho_{\bf B}}{2} \, \|{\bf B}_{n_B} - \tilde{\bf B}^k_{n_B} \|_F^2 \right) \cr
& \quad \qquad \left. +  \sum_{n_C=1}^{N_C} \left( {\bf Y}^k_{{\bf C}_{n_C}} * {\bf C}_{n_C} + \frac{\rho_{\bf C}}{2} \, \|{\bf C}_{n_C} - \tilde{\bf C}^k_{n_C} \|_F^2 \right)
\right)\cr
& \tilde{\bf A}_{n_A}^{k+1}  = \left( {\bf A}_{n_A}^{k+1} + \frac{1}{\rho_{\bf A}} {\bf Y}_{{\bf A}_{n_A}}^k \right)_+, ~n_A=1,\ldots,N_A, \cr
& \tilde{\bf B}_{n_B}^{k+1}  = \left( {\bf B}_{n_B}^{k+1} + \frac{1}{\rho_{\bf B}} {\bf Y}_{{\bf B}_{n_B}}^k \right)_+, ~n_B=1,\ldots,N_B,  \cr
&\tilde{\bf C}_{n_C}^{k+1}  = \left( {\bf C}_{n_C}^{k+1}  + \frac{1}{\rho_{\bf C}} {\bf Y}_{{\bf C}_{n_C}}^k \right)_+, ~n_C=1,\ldots,N_C, \cr
& {\bf Y}_{{\bf A}_{n_A}}^{k+1}  = {\bf Y}_{{\bf A}_{n_A}}^{k} + \rho_{\bf A} \left( {\bf A}_{n_A}^{k+1} - \tilde{\bf A}_{n_A}^{k+1} \right),
~n_A=1,\ldots,N_A, \cr
& {\bf Y}_{{\bf B}_{n_B}}^{k+1}  = {\bf Y}_{{\bf B}_{n_B}}^{k} + \rho_{\bf B} \left( {\bf B}_{n_B}^{k+1} - \tilde{\bf B}_{n_B}^{k+1} \right),
~n_B=1,\ldots,N_B, \cr
& {\bf Y}_{{\bf C}_{n_C}}^{k+1}  = {\bf Y}_{{\bf C}_{n_C}}^{k} + \rho_{\bf C} \left( {\bf C}_{n_C}^{k+1} - \tilde{\bf C}_{n_C}^{k+1} \right),
~ n_C=1,\ldots,N_C.
\end{split}
\label{ADMoM_NTF_vl}
\end{equation}
The minimization problem in the first line of (\ref{ADMoM_NTF_vl}) is non-convex.
Based on (\ref{f_X_vl}), we propose the alternating optimization scheme given in
(\ref{C_nC_sol}) at the next page.

\begin{table*}
\normalsize
\begin{equation}
\begin{split}
{\bf A}_{n_A}^{k+1} & = \underset{{\bf A}_{n_A}} {\rm argmin}
\left(
\left( \sum_{n_C=1}^{N_C} \frac{1}{2} \, \|{\bf X}^{(1)}_{n_A,n_C} - {\bf A}_{n_A}  ({\bf C}_{n_C}^k \odot {\bf B}^k)^T \|_F^2 \right)
+ {\bf Y}^k_{{\bf A}_{n_A}} * {\bf A}_{n_A} + \frac{\rho_{\bf A}}{2} \, \|{\bf A}_{n_A} - \tilde{\bf A}_{n_A}^k \|_F^2 \right) \cr
& = \left( \left( \sum_{n_C=1}^{N_C} {\bf X}^{(1)}_{n_A,n_C} ({\bf C}_{n_C}^k \odot {\bf B}^k) \right) +
\rho_{\bf A} \tilde{\bf A}_{n_A}^{k} - {\bf Y}_{{\bf A}_{n_A}}^{k} \right) \cr
& \qquad \qquad
\left( \left( \sum_{n_C=1}^{N_C} ({\bf C}_{n_C}^{k} \odot {\bf B}^k)^T ({\bf C}_{n_C}^{k} \odot {\bf B}^k) \right) +
\rho_{\bf A} {\bf I}_{F} \right)^{-1},
 ~\mbox{for}~n_A=1,\ldots,N_A, \cr
{\bf B}_{n_B}^{k+1} & = \underset{{\bf B}_{n_B}} {\rm argmin}
\left(
\left( \sum_{n_C=1}^{N_C} \frac{1}{2} \, \| {\bf X}^{(2)}_{n_B,n_C} - {\bf B}_{n_B}  ({\bf C}_{n_C}^k \odot {\bf A}^{k+1})^T \|_F^2 \right)
+ {\bf Y}^k_{{\bf B}_{n_B}} * {\bf B}_{n_B} + \frac{\rho_{\bf B}}{2} \, \|{\bf B}_{n_B}-\tilde{\bf B}_{n_B}^k \|_F^2 \right) \cr
& = \left( \left( \sum_{n_C=1}^{N_C} {\bf X}^{(2)}_{n_B,n_C} ({\bf C}_{n_C}^k \odot {\bf A}^{k+1}) \right)
+ \rho_{\bf B} \tilde{\bf B}_{n_B}^{k} - {\bf Y}_{{\bf B}_{n_B}}^{k}  \right)
\cr
& \qquad \qquad
\left( \left( \sum_{n_C=1}^{N_C} ({\bf C}_{n_C}^{k} \odot {\bf A}^{k+1})^T ({\bf C}_{n_C} ^{k} \odot {\bf A}^{k+1}) \right)
+ \rho_{\bf B} {\bf I}_{F} \right)^{-1},
~ \mbox{for}~  n_B=1,\ldots,N_B, \cr
{\bf C}_{n_C}^{k+1} & = \underset{{\bf C}_{n_C}} {\rm argmin}
\left(
\left( \sum_{n_B=1}^{N_B} \frac{1}{2} \, \|{\bf X}^{(3)}_{n_C,n_B} - {\bf C}_{n_C}  ({\bf B}_{n_B}^{k+1} \odot {\bf A}^{k+1})^T \|_F^2 \right)
+ {\bf Y}^k_{{\bf C}_{n_C}} * {\bf C}_{n_C} + \frac{\rho_{\bf C}}{2} \, \|{\bf C}_{n_C} - \tilde{\bf C}_{n_C}^k \|_F^2 \right) \cr
& = \left( \left( \sum_{n_B=1}^{N_B} {\bf X}^{(3)}_{n_C,n_B} ({\bf B}_{n_B}^{k+1} \odot {\bf A}^{k+1}) \right)
+ \rho_{\bf C} \tilde{\bf C}_{n_C}^{k} - {\bf Y}_{{\bf C}_{n_C}}^{k}  \right) \cr
& \qquad \qquad
\left( \left( \sum_{n_B=1}^{N_B} ({\bf B}_{n_B}^{k+1} \odot {\bf A}^{k+1})^T ({\bf B}_{n_B}^{k+1} \odot {\bf A}^{k+1}) \right)
+ \rho_{\bf C} {\bf I}_{F} \right)^{-1},
~ \mbox{for}~n_C=1,\ldots,N_C.
\end{split}
\label{C_nC_sol}
\end{equation}
\hrulefill
\end{table*}

Again, during each ADMoM iteration, we avoid the solution of constrained optimization problems.
Furthermore, and more importantly, having computed all algorithm quantities at iteration $k$,
the updates of ${\bf A}_{n_A}^k$, for $n_A=1,\ldots, N_A$, are {\em independent} and can be computed in {\em parallel}.
Then, we can compute in parallel the updates of ${\bf B}_{n_B}^k$, for $n_B=1,\ldots,N_B$,
and, finally, the updates of ${\bf C}_{n_C}^k$, for $n_C=1,\ldots,N_C$.

We note that we can solve problem (\ref{Problem_partitioned_NTF_ADMM_vl}) using the centralized ADMoM of Section
\ref{Section_ADMoM_NTF}. In fact, if we initialize the corresponding quantities of the two algorithms with
the same values, then the two algorithms evolve in {\em exactly} the same way. As a result, the study (for example, convergence analysis
and/or numerical behavior) of one of them is sufficient for the characterization of both.

Thus, via the distributed ADMoM, we simply uncover the {\em inherent parallelism}\/ in the updates of
the blocks of ${\bf A}^{k}$, ${\bf B}^{k}$, and ${\bf C}^{k}$.
In Appendix \ref{appendix_equivalence}, we present a detailed proof of the equivalence of these two forms of ADMoM NTF.

\subsection{A parallel implementation of ADMoM for large NTF}
\label{Subsection_architecture}

In the sequel, we briefly describe a simple implementation of ADMoM for large NTF on a mesh-type architecture.
In order to keep the presentation simple, we assume that (1) $N_A=N_B=N_C=N$ and
(2) each of the matrix unfoldings ${\bf X}^{(1)}$, ${\bf X}^{(2)}$, and ${\bf X}^{(3)}$ has been split
into $N^2$ blocks, with their $(i,j)$-th blocks stored at the $(i,j)$-th processing element,
for $i,j=1,\ldots,N$ (for related results in the matrix factorization context see
\cite{Gemulla:2011:LMF:2020408.2020426}). 

\begin{figure*}[t]
\centerline{\input{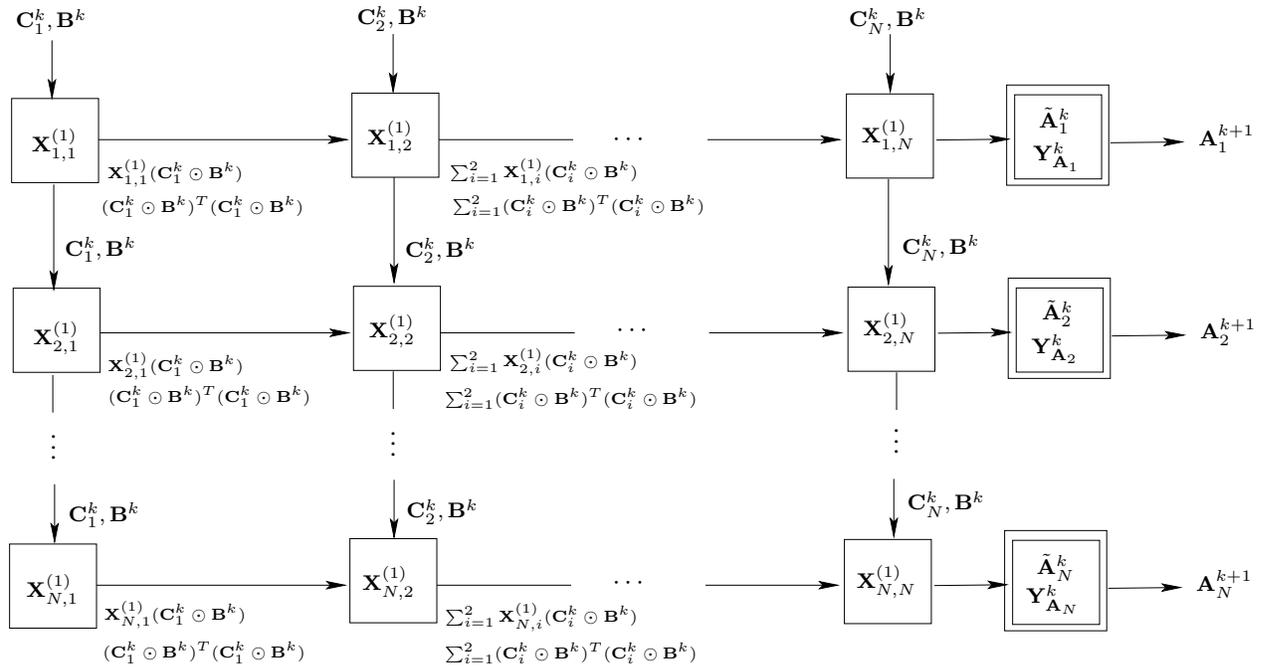}}
\caption{Distributed computation of ${\bf A}_{n}^{k+1}$, for $n=1,\ldots,N$.}
\label{fig_A_kp1}
\end{figure*}

In Figure \ref{fig_A_kp1}, we depict the data flow for the computation of the blocks of ${\bf A}^{k+1}$.
The inputs to the $N$ top processing elements are ${\bf C}_n^k$, for $n=1,\ldots,N$, as well as ${\bf B}^k$,
which is common input to all top processing elements. Each processing element uses its inputs and memory contents and
computes certain partial matrix sums. The communications between the processing elements are local and involve either the
forwarding of the terms ${\bf C}_n^k$, for $n=1,\ldots,N$, and ${\bf B}^k$ (top-down communication),
or the forwarding of the partial sums $\sum_{l=1}^j {\bf X}^{(1)}_{n,l} ({\bf C}_l^k \odot {\bf B}^k)$ and
$\sum_{l=1}^j ({\bf C}_l^k \odot {\bf B}^k)^T ({\bf C}_l^k \odot {\bf B}^k)$ (left-right communication), of dimensions
$\frac{I}{N} \times F$ and $F \times F$, respectively.
The computation of ${\bf A}_n^{k+1}$, for $n=1,\ldots,N$, amounts to solution of
$\rho$ systems of linear equations with common coefficient matrix
and takes place at the rightmost computing elements.

Then, using a similar strategy, we can compute the blocks of ${\bf B}^{k+1}$ and, finally, the
blocks of ${\bf C}^{k+1}$.
The updates of the auxiliary and dual variables are very simple and can be performed locally (see at the rightmost
computing elements of Figure \ref{fig_A_kp1}).

As we see in Figure \ref{fig_A_kp1}, in order to compute the blocks of the ${\bf A}^{k+1}$,
we use the appropriate blocks of ${\bf C}^k_n$, for $n=1,\ldots,N$, as well as the
{\em whole}\/ matrix ${\bf B}^{k}$. When the size of ${\bf B}^{k}$ is not very large, the
communication cost is {\em not}\/ prohibitive (analogous arguments holds
for the computation of the blocks of ${\bf B}^{k+1}$ and ${\bf C}^{k+1}$). Of course, if
one or more latent factors are very large, the communication cost significantly increases.

Concerning the distributed implementation of ADMoM for large tensor factorization with structural constraints
other than the non-negativity of the latent factors, we note that, if the size of
the latent factors is not very large, the communication cost of gathering together the blocks of
the auxiliary variables $\tilde{\bf A}^k$, $\tilde{\bf B}^k$, and $\tilde{\bf C}^k$,
is not prohibitive, enabling the computation of more complicated non-separable projections, like,
for example, projection onto the set of non-negative matrices with a certain maximum number of
non-negative elements.

Actual implementation of the distributed ADMoM for large NTF will depend on the specific parallel
architecture and programming environment used.
Since our aim in this paper is to introduce the basic methodology and computational framework, we leave those
customizations and performance tune-ups, which are further away from the signal processing core, for follow-up work
to be reported in the high-performance computing literature.

\section{Numerical Experiments}
\label{Section_Numerical}

\subsection{Comparison of ADMoM with NALS and NLS}
\label{subsection_comparison_with_NALS_NLS}

In our numerical experiments, we compare ADMoM NTF with (1) NALS NTF, as implemented in the ${\tt parafac}$
routine of the {\em N-way toolbox} for Matlab \cite{Nwaytoolbox_cite} and (2) NTF
using the nonlinear least-squares solvers (NLS), as implemented in the
${\tt sdf\_nls}$ routine of tensorlab \cite{Tensorlab_cite} (with the non-negativity option turned on in both cases).
In all cases, we use random initialization.
More specifically, the initialization of the ADMoM NTF is as follows. We give
non-negative {\em random}\/ values to ${\bf B}^0$ and ${\bf C}^0$ and zero values to 
the other state variables of the algorithm, namely, $\tilde{\bf A}^0$, $\tilde{\bf B}^0$,
$\tilde{\bf C}^0$, ${\bf Y}_{\bf A}^0$, ${\bf Y}_{\bf B}^0$, and 
${\bf Y}_{\bf C}^0$.\footnote{In certain cases, it may be possible to employ algebraic initialization schemes 
(e.g., see \cite{tomasi2006comparison} 
and references therein), but for NTF we observed that these perform (very) well only in (very) high SNR cases.}

\begin{table*}
\caption{Mean relative factorization error and mean and standard deviation of ${\tt cputime}$, 
in sec, for NALS, NLS and ADMoM NTF.
}
\centering
\begin{tabular}{crrrrrrrrr}
\hline
\mbox{Size} & $F$ & $\sigma_N^2$~ & $\begin{array}{c} {\rm mean({\tt RFE})} \end{array}$
& $\begin{array}{c} {\tt NALS} \\ {\rm mean({\tt t})} \end{array}$
& $\begin{array}{c} {\tt NLS} \\ {\rm mean({\tt t})} \end{array}$ &
$\begin{array}{c} {\tt ADMoM} \\ {\tt mean(t)} \end{array}$ & $\begin{array}{c} {\tt NALS} \\ {\tt std(t)} \end{array}$
& $\begin{array}{c} {\tt NLS} \\ {\rm std({\tt t})} \end{array}$ &
$\begin{array}{c} {\tt ADMoM} \\ {\tt std(t)} \end{array}$  \\
\hline
$3000 \times 50 \times 50$ & $3$  & $10^{-2}$ & $0.2156$~ & $12.4010$~  & $17.2584$~  & $7.1806$~   & $1.4770$~ & $6.9348$~ & $3.8647$~   \\
                           &      & $10^{-4}$ & $0.0221$~ & $16.6500$~  & $16.5962$~  & $7.5098$~   & $1.9625$~ & $3.3589$~ & $4.2499$~ \\
                           & $30$ & $10^{-2}$ & $0.0260$~ & $212.0598$~ & $115.7030$~  & $110.4128$~ & $11.4579$~ & $8.0409$~ & $91.4882$~  \\
                           &      & $10^{-4}$ & $0.0026$~  & $270.6674$~ & $117.2152$~ & $148.5052$~ & $11.5324$~ & $11.3154$~ & $149.0542$~   \\
\hline
$400 \times 400 \times 50$ & $3$  & $10^{-2}$ & $0.2175$~ & $8.6174$~  & $4.7124$~ & $8.0710$~ & $0.9264$~ & $1.4596$~ & $4.4952$~ \\
                           &      & $10^{-4}$ & $0.0222$~ & $11.0670$~ & $4.8916$~ & $8.3614$~ & $1.2075$~ & $1.6443$~ & $3.1833$~  \\
                           & $30$ & $10^{-2}$ & $0.0260$~ & $71.1950$~ & $31.7546$~ & $106.1362$~ & $8.2119$~ & $3.2162$~ & $76.5743$~ \\
                           &      & $10^{-4}$ & $0.0026$~  & $92.4838$~ & $31.1690$~ & $94.5734$~  & $7.9280$~ & $3.7062$~ & $88.1671$~ \\
\hline
$200 \times 200 \times 200$ & $5$  & $10^{-2}$ & $ 0.1400$~ & $10.9142$~ & $4.1756$~ & $12.5190$~  & $1.2783$~ & $0.6915$~ & $2.1334$~ \\
                            &      & $10^{-4}$ & $0.0143$~  & $14.2882$~  & $4.1070$~ & $12.6184$~  & $2.5817$~ & $0.9059$~ & $2.1616$~ \\
                            & $30$ & $10^{-2}$ & $0.0260$~  & $55.0806$~  & $16.1838$~ & $33.9268$~ & $4.4886$~ & $1.3649$~ & $12.2687$~ \\
                            &      & $10^{-4}$ & $0.0026$~  & $70.0238$~  & $16.7624$~ & $32.1670$~ & $6.5737$~ &$1.4152$~ & $10.0182$~  \\
\hline
\end{tabular}
\label{Table_1}
\end{table*}

In extensive numerical experiments, we have observed that the relative performance of the algorithms depends on
the size and rank of the tensor as well as the additive noise power. Thus,
we consider $12$ different scenarios, corresponding to the combinations of the following cases:
\begin{enumerate}
\item one, two, or three tensor dimensions are large;
\item rank $F$ is small or large;
\item additive noise is weak or strong.
\end{enumerate}
For each scenario, we generate $R=50$ realizations of tensor $\underline{\bf X}$ as follows.
We generate random matrices ${\bf A}^o$, ${\bf B}^o$, and ${\bf C}^o$ with
i.i.d. ${\cal U}[0,1]$ elements (using the ${\tt rand}$ command of Matlab)
and construct
$\underline{\bf X} = [{\bf A}^o, {\bf B}^o, {\bf C}^o] + \underline{\bf N}$,
where $\underline{\bf N}$ consists of i.i.d.  ${\cal N}(0,\sigma_N^2)$  elements.
For each realization, we solve the NTF problem with (1) NALS (${\tt parafac}$), (2) NLS (${\tt sdf\_nls}$), 
and (3) ADMoM.

We designed our experiments so that, upon convergence,  all algorithms achieve practically
the same relative factorization error. Towards this end, we set the values of the stopping parameters
as follows: the parameter ${\tt Options(1)}$ of ${\tt parafac}$ is set
to ${\tt Options(1)}=10^{-5}$, the parameter ${\tt TolFun}$ of ${\tt sdf\_nls}$ is set to ${\tt TolFun}=10^{-8}$, and
the ADMoM stopping parameters are set to $\epsilon^{\rm abs}=10^{-4}$ and $\epsilon^{\rm rel}=10^{-4}$.

In all cases, the initial values of the ADMoM penalty terms are $\rho_{\bf M}=1$,
for ${\bf M}={\bf A}, {\bf B}, {\bf C}$, while the ADMoM penalty term adaptation parameters are $\mu = 8$,
$\tau^{\rm incr}=4$, $\tau^{\rm decr}=2$.

In practice, convergence properties of ADMoM NTF depend on the (random) initialization point.
In some cases, convergence may be quite fast while, in others, it may be quite slow. As we shall see in the sequel,
this phenomenon seems more prominent in the cases where rank $F$ is large. In order to overcome the slow
convergence properties associated
with bad initial points, we adopted the following strategy.
We execute ADMoM NTF for up to $n_{\max}=400$ iterations (we have observed that,
in the great majority of the cases in the scenarios we examined,
this number of iterations is sufficient for convergence when we start from a
good initial point). If ADMoM does not converge within $n_{\max}$
iterations, then we restart it from another random initial point; we repeat this procedure until ADMoM converges.\footnote{Of
course, one may think of more elaborate strategies  such as, for example, running in parallel more
than one versions of the algorithm, with different initializations.}

\begin{figure}
\centerline{\epsfsize=0.6\hsize\epsfbox{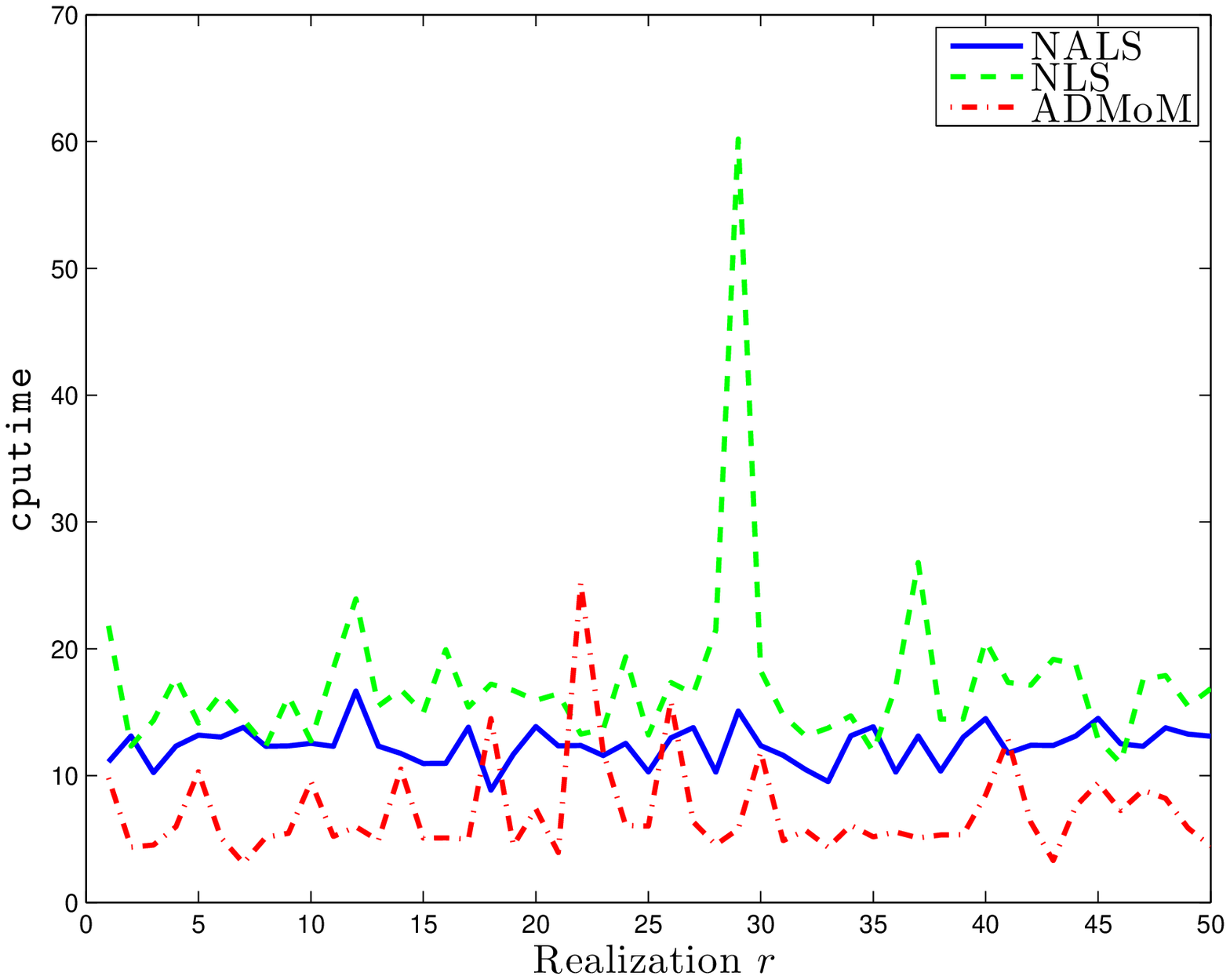}}
\caption{${\tt cputime}$ for $I=3000$\, $J=K=50$, $F=3$, and $\sigma_N^2=10^{-2}$.
NALS (blue solid line), NLS (green dashed line), ADMoM (red dotted-dashed line).}
\label{Fig_cputime_3000_50_50_3_1}
\end{figure}

\begin{figure}
\centerline{\epsfsize=0.6\hsize\epsfbox{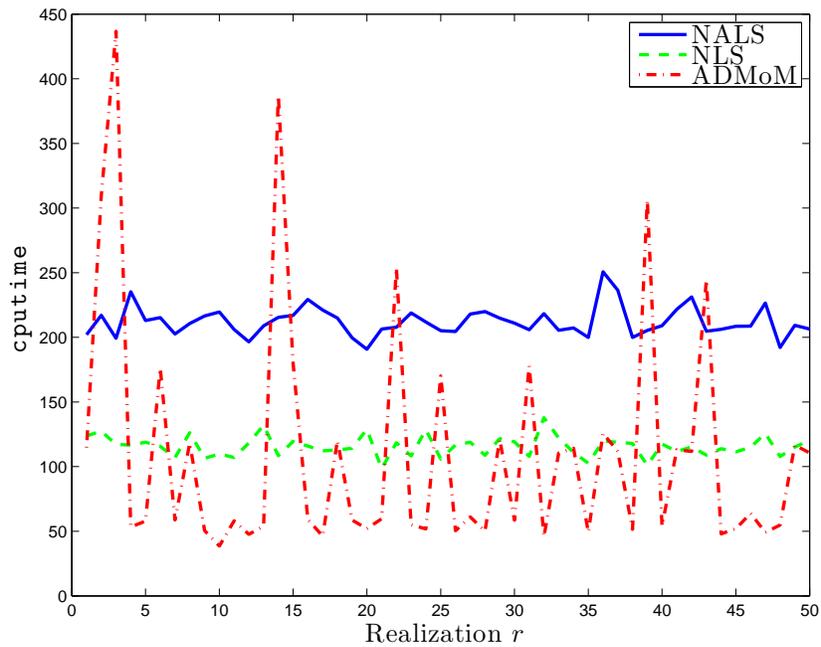}}
\caption{${\tt cputime}$ for $I=3000$, $J=K=50$, $F=30$ and $\sigma_N^2=10^{-2}$.
NALS (blue solid line), NLS (green dashed line), ADMoM (red dotted-dashed line).}
\label{Fig_cputime_3000_50_50_30_1}
\end{figure}

Before proceeding, we mention that {\em all}\/  the 
algorithms converged in {\em all}\/ the realizations we run.

Since an accurate statement about the computational complexity per iteration
of ${\tt parafac}$ is not easy, the metric we used for comparison of the algorithms is the ${\tt cputime}$
of Matlab. Despite the fact that ${\tt cputime}$ is strongly  dependent
on the computer hardware and the actual algorithm implementation, we feel that it is a useful metric
for the assessment of the relative efficiency of the algorithms.\footnote{For our experiments, we
run Matlab 2014a on a MacBook Pro with a $2.5$ GHz Intel Core i7 Intel processor and $16$ GB RAM.} 
The reason is that we used carefully developed, publicly available Matlab toolbox implementations of the 
baseline algorithms, and we carefully coded our ADMoM NTF implementation. 

In Table \ref{Table_1}, we present the mean and standard deviation of ${\tt cputime}$, in seconds,
denoted as ${\tt mean(t)}$ and ${\tt std(t)}$, respectively, for NALS, NLS, and ADMoM.
We also present the mean relative factorization error (which is common to all algorithms up to four decimal digits), 
defined as 
\[
{\tt mean(RFE)}:= \frac{1}{n}\sum_{k=1}^n \frac{\|\underline{\bf X}_k -
[{\bf A}_k, {\bf B}_k, {\bf C}_k\|_F}{\|\underline{\bf X}_k\|_F},
\]
where $\underline{\bf X}_k$ is the $k$-th noisy tensor realization 
and ${\bf A}_k$, ${\bf B}_k$, and ${\bf C}_k$ are the factors returned by a
factorization algorithm.
Our observations are as follows:
\begin{enumerate}
\item
There is no clear winner. Certainly, for high ranks, NLS has very good behavior.

\item
In general, both NALS and NLS have more predictable behavior than ADMoM. Especially for high ranks,
the ${\tt cputime}$ of our implementation of ADMoM has large variance.

\item
For small ranks, ADMoM looks more competitive and, in the cases where one dimension
is much larger than the other two, it behaves very well (we shall say more on this later).
\end{enumerate}
In order to get a better feeling of the behavior of the three algorithms,
we plot their ${\tt cputime}$, along the $50$ realizations we used to
obtain the averages of Table \ref{Table_1}, for two different scenarios.
In Figure \ref{Fig_cputime_3000_50_50_3_1}, we
consider the case for $I=3000$, $J=K=50$, $F=3$ and $\sigma_N^2=10^{-2}$. We
observe that the behavior of the algorithms is stable, in the sense that there is a clear ordering among
the three algorithms, with no large variations.
In Figure \ref{Fig_cputime_3000_50_50_30_1}, we keep the dimensions and the noise power the same as before and
increase the rank to $F=30$. We observe that the variance of ADMoM ${\tt cputime}$ has significantly increased,
while both NALS and NLS show stable behavior.
When ADMoM starts from a good initial point, it converges faster than NALS and NLS while, when it starts 
from  bad initial points, it needs one or more restarts.

\begin{table*}
\caption{Mean relative factorization error and mean and standard deviation of ${\tt cputime}$, in sec, for NALS, NLS and ADMoM NTF for
$I=10^4$, $J=K=50$, $F=10$, and $\sigma_N^2=10^{-2}$.}
\centering
\begin{tabular}{ccccccccc}
\hline
$\begin{array}{c} {\tt NALS} \\ {\tt mean(RFE)} \end{array}$ &
$\begin{array}{c} {\tt NLS} \\ {\tt mean(RFE)} \end{array}$ &
$\begin{array}{c} {\tt ADMoM} \\ {\tt mean(RFE)} \end{array}$ &
$\begin{array}{c} {\tt NALS} \\ {\tt mean(t)} \end{array}$ &
$\begin{array}{c} {\tt NLS} \\ {\tt mean(t)} \end{array}$ &
$\begin{array}{c} {\tt ADMoM} \\ {\tt mean(t)} \end{array}$ &
$\begin{array}{c} {\tt NALS} \\ {\tt std(t)} \end{array}$ &
$\begin{array}{c} {\tt NLS} \\ {\tt std(t)} \end{array}$  &
$\begin{array}{c} {\tt ADMoM} \\ {\tt std(t)} \end{array}$  \\
\hline
$0.0752$ & $0.0758$  & $0.0751$ & $116.0722$& $201.6556$  & $30.6460$  & $19.8483$ & $38.5401$ & $16.5823$ \\
\hline
\end{tabular}
\label{Table_2}
\end{table*}

\begin{figure}[t]
\centerline{\epsfsize=0.6\hsize\epsfbox{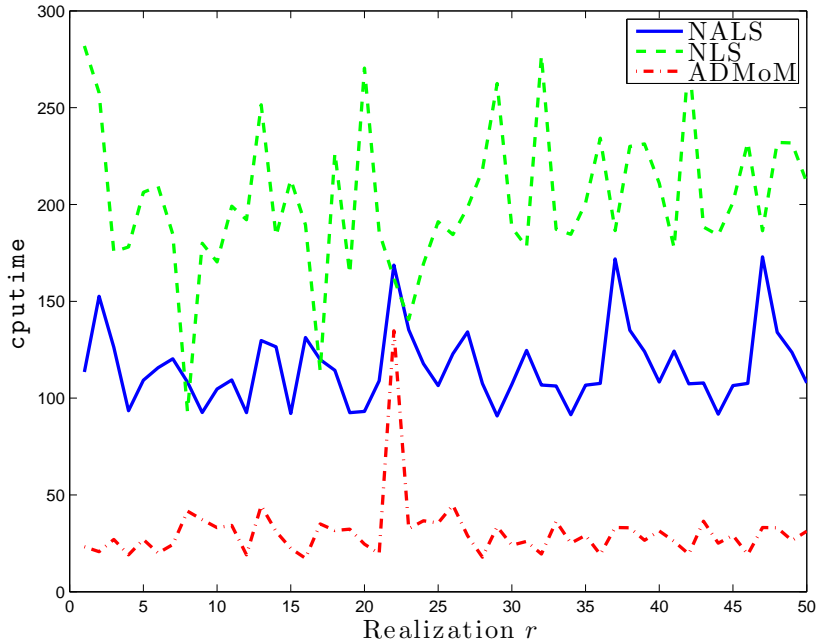}}
\caption{${\tt cputime}$ for $I=10^4$, $J=K=50$, $F=10$, and $\sigma_N^2=10^{-2}$.
NALS (blue solid line), NLS (green dashed line), ADMoM (red dotted-dashed line)..}
\label{Fig_cputime_10000_50_50_10_1}
\end{figure}

In order to check if ADMoM maintains its advantage over NALS and NLS in the cases where one dimension is
very large, compared with the other two, and the rank is relatively small, we performed an experiment
with $I=10^4$, $J=K=50$, $F=10$, and $\sigma_N^2=10^{-2}$. However, in this case, we used somewhat relaxed
stopping conditions for all algorithms; more specifically, we used 
${\tt Options(1)}=10^{-3}$, ${\tt TolFun}=10^{-6}$,
$\epsilon^{\rm abs}=10^{-3}$, and $\epsilon^{\rm rel}=10^{-3}$. In Table \ref{Table_2}, we present the
mean relative factorization errors and the mean and standard deviation of ${\tt cputime}$.
As we can see, both NALS and NLS are slightly less accurate than ADMoM, in terms of relative factorization error,
which means that their stopping criteria are more relaxed. In terms of ${\tt cputime}$, we see that ADMoM is much faster than both
NALS and NLS. In Figure \ref{Fig_cputime_10000_50_50_10_1}, we plot the ${\tt cputime}$ of the three
algorithms for the $50$ realizations of the experiment. Again, we see the significant difference between ADMoM and
both NALS and NLS. We note that if we had used as values of the stopping parameters those of our initial experiments, then the gain of
ADMoM, compared with NALS and NLS, would have been much greater. However, we believe that we have made clear that, in this case,
ADMoM has a clear advantage. We have made analogous observations for larger $I$.

\subsection{A closer look at ADMoM}
\label{subsection_closer_look_ADMoM}

\begin{figure}[t]
\centerline{\epsfsize=0.6\hsize\epsfbox{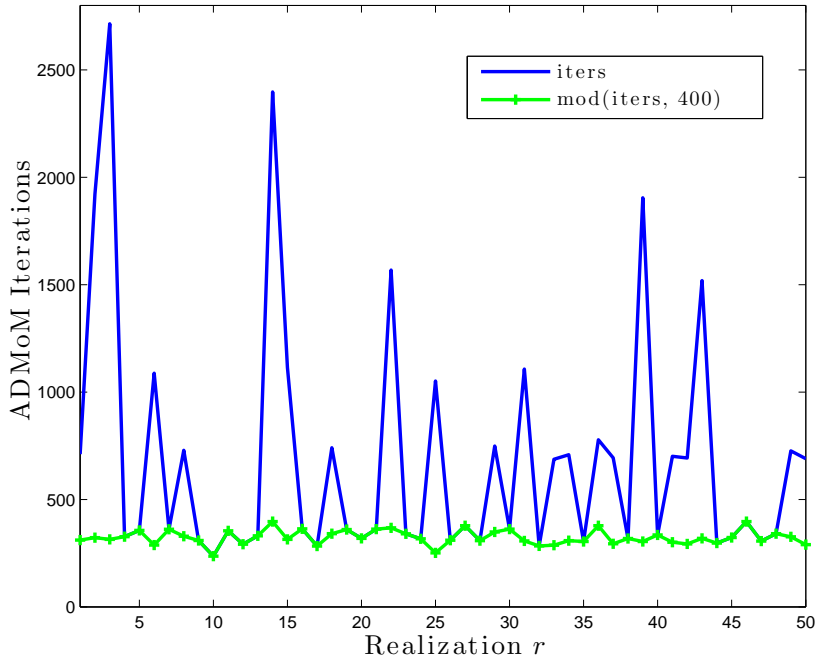}}
\caption{Number of ADMoM iterations for $I=3000$, $J=50=K=50$, $F=30$, and $\sigma_N^2=10^{-2}$;
${\tt iters}$ (blue line), and  ${\rm mod}({\tt iters}, n_{\max})$ (green line).}
\label{Fig_ADMoM_iters_3000_50_50_30_1}
\end{figure}

In order to get a more detailed view of the convergence properties of ADMoM, we return to the scenario with
$I=3000$, $J=K=50$, $F=30$, and $\sigma_N^2=10^{-2}$, whose ${\tt cputime}$ we plot
in Figure \ref{Fig_cputime_3000_50_50_30_1}. We recall that, in order to converge in this case, ADMoM needed often restarts.
In Figure \ref{Fig_ADMoM_iters_3000_50_50_30_1},
we plot the total number of ADMoM iterations, denoted as ${\tt iters}$, and the number of ADMoM iterations during
its final way to convergence, which is equal to
${\rm mod}({\tt iters}, n_{\max})$. As expected, ${\tt iters}$ is compatible with the corresponding ${\tt cputime}$
(see the red line in Figure \ref{Fig_cputime_3000_50_50_30_1}).
Quantity ${\rm mod}({\tt iters}, n_{\max})$ shows how many iterations are required for convergence
 if ADMoM always starts from good initial points. We observe that ${\rm mod}({\tt iters}, n_{\max})$ is quite
 stable around its mean, which is approximately equal to $320$. This gives an estimate of the fastest possible
 ADMoM convergence in this case.

\subsection{ADMoM NTF with under- and over-estimated rank}

\begin{figure}[t]
\centerline{\epsfsize=0.6\hsize\epsfbox{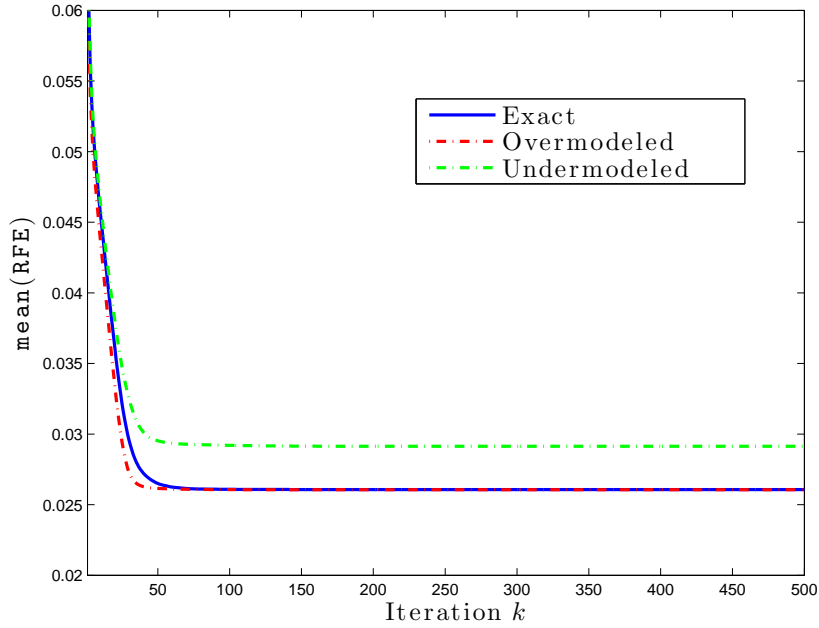}}
\caption{Average relative factorization errors for $I=J=K=100$, $F=30$, and $\sigma_N^2=0.1$.
Exact rank case (solid blue line), over-estimated rank by $1$ (dotted-dashed red line),
underestimated rank by $1$ (green dashed line).}
\label{Fig_noisy_rel_fact_err}
\end{figure}

\begin{figure}[t]
\centerline{\epsfsize=0.6\hsize\epsfbox{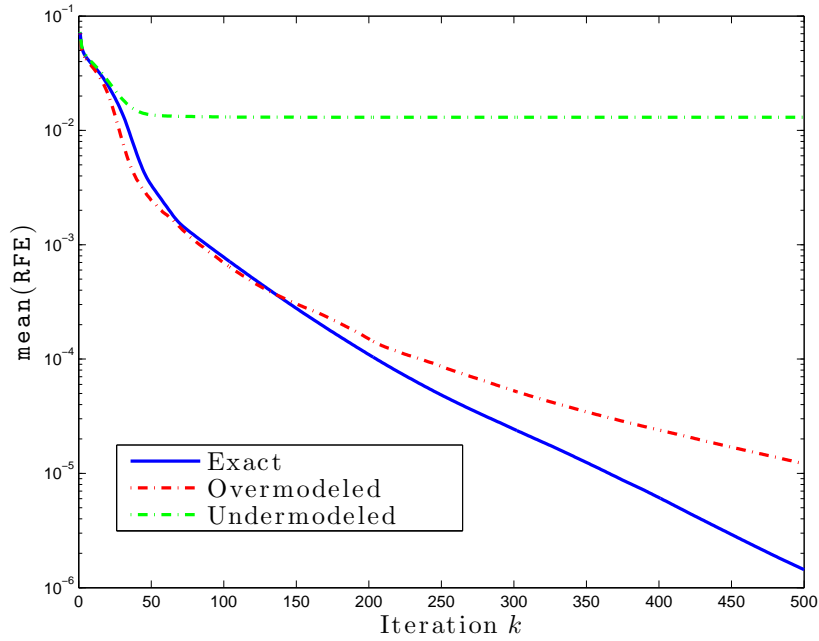}}
\caption{Average relative factorization errors for  noiseless case with $I=J=K=100$ and $F=30$.
Exact rank case (solid blue line), over-estimated rank by $1$ (dotted-dashed red line),
underestimated rank by $1$ (green dashed line).}
\label{Fig_noiseless_rel_fact_err}
\end{figure}

In the sequel, we consider ADMoM behavior in the cases where we under- or over-estimate the true rank, in both noisy and
noiseless cases. Towards this end, we fix $I=J=K=100$ and $F=30$, and investigate ADMoM with exact rank
as well as with rank under- and over-estimated by $1$. We expect that, in this case, all versions of ADMoM may need restarts.
In the sequel, we examine the influence of under- and over-estimating the rank on (1) factorization accuracy and (2)
number of restarts. Of course, in under-modeled cases, we expect that the relative factorization error will be higher than that of the
true rank case. However, we know nothing in advance about ADMoM behavior in over-modeled cases. In order to get insight into these issues,
we perform the following experiment. We set stopping parameters
$\epsilon^{\rm abs}=\epsilon^{\rm rel}=10^{-4}$ and run each of the three versions of ADMoM for
$n_{\rm max}=500$ iterations. For each version, we proceed as follows: if it converges
within $n_{\rm max}$ iterations, we stop; otherwise, we restart, and repeat until convergence. Thus, finally, the
number of iterations for each ADMoM version will be a multiple of $n_{\rm max}$. For the computation of
the trajectory of the mean relative factorization error we use only the last $n_{\rm max}$ values; in this way,
we avoid the influence of bad initial points. However, we keep count of the restarts of each version and, thus,
can assess the time it needs to achieve convergence.

In Figure \ref{Fig_noisy_rel_fact_err}, we plot the average relative factorization
errors (computed over $50$ realizations in the way we mentioned before), versus the iteration number, for $\sigma_N^2 = 10^{-2}$.
As was expected, the ADMoM version with under-estimated rank converges to a higher relative factorization error.
We observe that the average relative factorization errors for ADMoM with exact rank and rank over-estimated by $1$
follow almost the same trajectory. The average numbers of restarts for the three ADMoM versions are
${\tt mean}({\rm restarts}_{\tt exact})=1.08$, ${\tt mean}({\rm restarts}_{\tt under})=1.22$, and
${\tt mean}({\rm restarts}_{\tt over})=1.92$. Thus, in the cases of over-estimated rank, we finally achieve
a relative factorization error trajectory as good as in the exact rank case, but we may need more restarts and,
thus, more time. This implies that the probability of bad initial  points may increase.

In Figure \ref{Fig_noiseless_rel_fact_err}, we plot the same quantities for noiseless data.
Again, the ADMoM behavior in the under-modeled case is as expected. Interestingly, we observe that
there is no relative factorization error floor neither for the exact rank nor for the over-estimated by $1$ rank case.
Reasonably, after a certain precision level, the over-modeled case converges slower. The average numbers of
restarts for the three ADMoM versions are ${\tt mean}({\rm restarts}_{\tt exact})=1.14$,
${\tt mean}({\rm restarts}_{\tt under})=1.54$, and ${\tt mean}({\rm restarts}_{\tt over})=1.06$.

We have observed similar behavior for more drastic rank under- and over-estimation.

\subsection{ADMoM for NTF with box-linear constraints}

\begin{figure}[t]
\centerline{\epsfsize=0.6\hsize\epsfbox{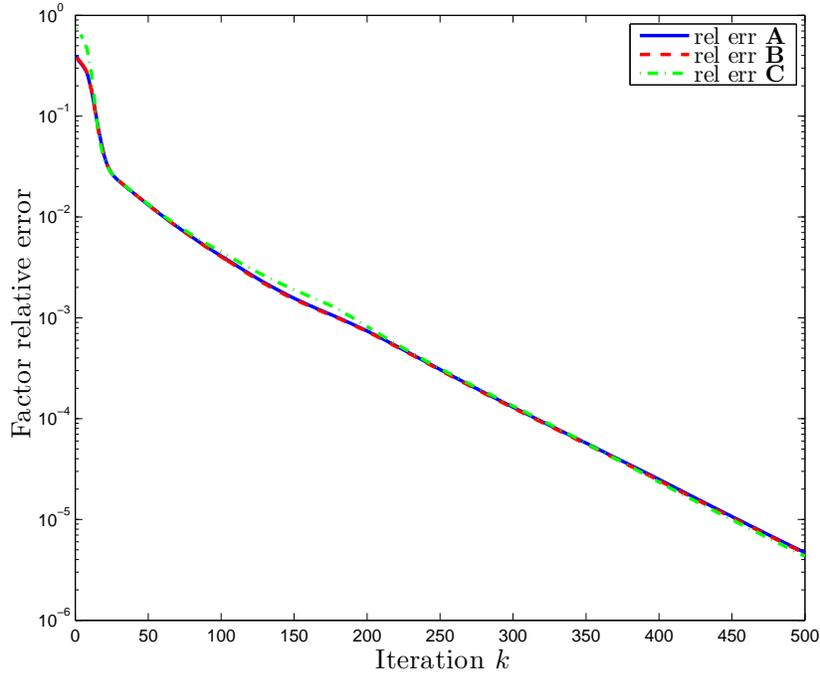}}
\caption{Average relative factorization errors for  noiseless case with $I=J=100$, $K=50$ and $F=20$.
Factors ${\bf A}$ (solid blue line), ${\bf B}$  (dashed red line), and ${\bf C}$ (green dotted-dashed line).}
\label{Fig_noiseless_Box_Linear}
\end{figure}

In our final experiment, we briefly consider NTF for the case where two of the latent factors, 
say ${\bf A}$ and ${\bf B}$, are non-negative while ${\bf C}$ is subject to box-linear constraints in the 
sense that each row of ${\bf C}$ is a probability mass function, that is, has non-negative elements with sum 
equal to 1.

The only difference between the ADMoM for this case and the ADMoM for NTF is that,
instead of computing ${\bf C}^{k+1}$ as the solution of an unconstrained least-squares problem, we compute 
it as the solution of linearly constrained least-squares; note that both cases exhibit closed-form solutions. 

In Figure \ref{Fig_noiseless_Box_Linear}, we illustrate the behavior of ADMoM in this case 
by plotting the trajectories of the norms of the average (over $50$ realizations) relative estimation errors
of the latent factors, as computed by function ${\tt cpderr}$ of tensorlab, versus the iteration number, for a noiseless
case with $I=J=100$, $K=50$ and $F=5$. We observe that ADMoM works to very high precision.

\begin{table*}[t]
\normalsize
\begin{equation}
\begin{split}
\left[ \begin{array}{c} {\bf A}^{k+1}_1 \cr \vdots \cr {\bf A}^{k+1}_{N_A}  \end{array} \right]
& =  \left(
\left[ \begin{array}{ccc} {\bf X}^{(1)}_{1,1} & \cdots & {\bf X}^{(1)}_{1,N_C} \cr
\vdots & \ddots & \vdots \cr
{\bf X}^{(1)}_{N_A,1} & \cdots & {\bf X}^{(1)}_{N_A,N_C}
\end{array} \right]
\left[ \begin{array}{c} {\bf C}^k_1 \odot {\bf B}^k \cr \vdots \cr {\bf C}^k_{N_C} \odot {\bf B}^k  \end{array} \right] +
\left[ \begin{array}{c} \rho_{\bf A} \tilde{\bf A}^k_1 - {\bf Y}^k_{{\bf A}_1} \cr \vdots \cr
\rho_{\bf A} \tilde {\bf A}^{k}_{N_A} - {\bf Y}^k_{{\bf A}_{N_A}}  \end{array} \right]
\right)
\cr
& \qquad \qquad \qquad \times
\left( \left( \sum_{n_C=1}^{N_C} ({\bf C}_{n_C}^k\odot {\bf B}^k )^T  ({\bf C}_{n_C}^k\odot {\bf B}^k ) \right) + \rho_{\bf A} {\bf I}_F \right)^{-1}.
\end{split}
\label{partitioned_update_Ak}
\end{equation}
\hrule
\end{table*}

\subsection{Discussion}
\label{subsection_discussion}

Our numerical results are encouraging and suggest that, in many cases,
ADMoM NTF can efficiently achieve close to state-of-the-art factorization accuracy.
The fact that ADMoM  is suitable for high-performance parallel implementation (the first NTF algorithm with
this property, as far as we know) can only increase its potential. Thus, we believe that it
will be a valuable tool in the NTF toolbox.

Obviously, in order to fully uncover the pros and cons of ADMoM NTF, more extensive experimentation is required.
But our intention in this paper is to give the fundamental ideas and some basic performance metrics. Experiments
with real-world data (using ADMoM for tensor completion and factorization) as well as constraints well beyond 
non-negativity are ongoing work.

A weak point of the version of ADMoM we developed in this manuscript is the high ${\tt cputime}$ variance
in cases of high rank. The improvement of the behavior of ADMoM in these cases remains
a very interesting problem.  To achieve this goal, it might be possible to combine
elements of NLS and ADMoM and derive a more efficient algorithm. However, more research efforts are needed
in this direction.

As we mentioned, if the centralized and the distributed algorithms start from the same
initial point, they evolve in exactly the same way. Thus, distributed ADMoM inherits
the convergence properties of centralized ADMoM.

\section{Conclusion}
\label{Section_Conclusion}

Motivated by emerging big data applications, involving multi-way tensor data, and the ensuing need for
scalable tensor factorization tools, we developed a new constrained tensor factorization framework based
on the ADMoM. We used non-negative factorization of third order tensors as an example to work out the main
ideas, but our approach can be generalized to higher order tensors, many other types of constraints on
the latent factors, as well as other tensor factorizations and tensor completion.
Our numerical experiments were encouraging, indicating that, in many cases, the ADMoM-based NTF has
high potential as an alternative to the state-of-the-art and, in some cases, it may become state-of-the-art.
The fact that it is naturally amenable to parallel implementation can only increase its potential.
The improvement of its behavior in the high rank cases remains a very interesting problem.

\appendices

\section{Extension to higher order tensors}
\label{partitionings_extension}

In this appendix, we highlight how our approach can be extended to higher order tensors.
We focus on fourth-order tensors, with the general case being obvious.
If $\underline{\bf W}=[{\bf A}, {\bf B}, {\bf C}, {\bf D}]$, then its matrix unfoldings satisfy relations
\begin{equation*}
\begin{split}
{\bf W}^{(1)} & = {\bf A} \,({\bf D} \odot ({\bf C} \odot {\bf B})  )^T,\cr
{\bf W}^{(2)} & = {\bf B} \, ({\bf D} \odot ({\bf C} \odot {\bf A})  )^T, \cr
{\bf W}^{(3)} & = {\bf C} \, ({\bf D} \odot ({\bf B} \odot {\bf A})  )^T, \cr
{\bf W}^{(4)} & = {\bf D} \, ({\bf C} \odot ({\bf B} \odot {\bf A})  )^T.
\end{split}
\end{equation*}
Partitioning matrices ${\bf A}$, ${\bf B}$, ${\bf C}$, and ${\bf D}$ as in subsection \ref{subsection_partitioned_matrix_unfoldings},
we obtain that matrix ${\bf W}^{(1)}$ can be partitioned into $N_A \times N_D$ blocks, with the $(i,j)$-th block
being equal to
\begin{equation*}
{\bf W}^{(1)}_{i,j} = {\bf A}_i \,({\bf D}_j \odot ({\bf C} \odot {\bf B}))^T.
\end{equation*}
Analogous partitionings apply to the other matrix unfoldings. Then, development of ADMoM NTF
(centralized and distributed) is rather easy.

\section{On the equivalence of the centralized and the distributed ADMoM NTF}
\label{appendix_equivalence}

A simple proof of the equivalence of the centralized and the distributed ADMoM NTF is as follows.
We focus on the update of ${\bf A}^k$ of the centralized algorithm and the updates of its blocks, ${\bf A}^k_{n_A}$, for
$n_A=1,\ldots,N_A$, of the distributed algorithm, and prove that they are equivalent.
We remind that 
\begin{equation*}
\begin{split}
{\bf A}^{k+1} & = \left( {\bf X}^{(1)} ({\bf C}^k\odot {\bf B}^k) +
\rho_{\bf A} \tilde{\bf A}^k - {\bf Y}_{\bf A}^k \right) \cr
& \hspace{2cm}
\left(  ({\bf C}^k\odot {\bf B}^k)^T ({\bf C}^k\odot {\bf B}^k) + \rho_{\bf A} {\bf I}_F\right)^{-1}.
\end{split}
\end{equation*}
Using the partitionings of ${\bf C}^k\odot {\bf B}^k$ (see subsection \ref{subsection_partitioned_matrix_unfoldings}),
it can be shown that
\begin{equation*}
({\bf C}^k\odot {\bf B}^k)^T ({\bf C}^k\odot {\bf B}^k) =
\sum_{n_C=1}^{N_C} ({\bf C}_{n_C}^k\odot {\bf B}^k )^T  ({\bf C}_{n_C}^k\odot {\bf B}^k ).
\end{equation*}
Rewriting the update of ${\bf A}^k$ in terms of partitioned matrices, we obtain
(\ref{partitioned_update_Ak}) at the top of this page.
If we focus on a certain block of ${\bf A}^{k+1}$ in (\ref{partitioned_update_Ak}),
then we obtain the corresponding update of the distributed algorithm (see (\ref{C_nC_sol})).
We observe that the matrix inverse in the second line of (\ref{partitioned_update_Ak}) is {\em common}\/
to all blocks, and should be computed once.

Analogous statements hold for the updates of ${\bf B}^k$ and ${\bf C}^k$. The equivalence of
the updates of the rest of the variables is trivial.

Thus, in fact, using the partitionings of subsection \ref{subsection_partitioned_matrix_unfoldings}, the distributed ADMoM
simply uncovered the inherent parallelism of the centralized ADMoM.

\bibliographystyle{IEEEtran}
\bibliography{ADMoM_NTF_bib}

\end{document}